\input amstex
\documentstyle{amsppt}
\topmatter

\title Hamiltonian selfdistributive quasigroups\endtitle

\author Dolors Herbera, Tom\'a\v s Kepka, and Petr N\v emec\endauthor

\address Department de Matem\`atiques, Universitat Aut\`onoma de Bar\-ce\-lona,
08193 Be\-lla\-te\-rra (Barcelona), Spain
\endaddress
\email dolors\@mat.uab.es\endemail

\address Department of Algebra, MFF UK, Sokolovsk\'a
83, 186 75 Praha 8, Czech Republic
\endaddress
\email kepka\@karlin.mff.cuni.cz\endemail

\address Department of Mathematics, \v CZU, Kam\'yck\'a 129,
165 21 Praha 6 -- Suchdol, Czech Republic\endaddress
\email nemec\@tf.czu.cz\endemail

\thanks The research of the first author was partially supported
by the DGI
and the European Regional Development Fund, jointly, through
Project BFM2002-01390, and by the Comissionat per Universitats i
Recerca of the Generalitat de Ca\-ta\-lunya.
The second author was supported by the institutional
grant MSM 113200007 and by the Grant Agency of Charles University,
grant \#269/2001/B-MAT/MFF, the third author by
the Grant Agency of Czech Republic, grant \#GA\v CR-201/02/0594.
\endthanks



\abstract The problem of the existence of non-medial distributive
hamiltonian quasigroups is solved.
Translating this problem first to commutative Moufang loops with
operators, then to ternary algebras and, finally, to cocyclic
modules
over\linebreak $\Bbb Z[x,x^{-1},(1-x)^{-1}]$, it
is shown that every non-medial distributive hamiltonian
quasigroup has at least 729 elements and that there are just two
isomorphism classes of such quasigroups of the least cardinality.
The quasigroups representing these two classes are
anti-isomorphic.
\endabstract

\endtopmatter

\document
{
\define\bo{\bold}
\define\R{\bo R}
\define\F{\boldsymbol\Phi}
\define\I{\bo I}
\define\SS{\bo S}
\define\T{\bo T}
\redefine\P{\bo P}
\define\Hom{\roman{Hom}}
\define\Ker{\roman{Ker}}
\define\An{\roman{An\,}}
\define\id{\roman{id}}
\define\Z{\roman{Z}}
\define\A{\roman{A}}
\define\Cend{\roman{Cend}}
\define\gen{\roman{gen}}
\define\K{\roman{K}}
\redefine\S{\roman{S}}
\define\J{\roman{J}}
\redefine\*{\circledast}
\redefine\Im{\roman{Im}}
\define\x{\bo x}
\define\y{\bo y}
\define\q{\roman q}
\define\U{\bo u}
\define\V{\bo v}
\define\W{\bo w}
\define\und{\underline}
\define\si{\sigma}
\define\th{\vartheta}
\define\la{\lambda}
\define\al{\alpha}
\define\be{\beta}
\define\ga{\gamma}
\define\ro{\varrho}
\define\ka{\varkappa}

\head 0. Introduction\endhead

The first explicit allusion to the left and right instances of
selfdistributivity (i.e., $(x(yz)\bumpeq(xy)(xz)$ and
$(xy)z\bumpeq(xz)(yz)$) seems to appear in \cite{39} where one
can read the following comment: {\it``These are other cases of
the distributive principle. \dots These formulae, which have
hitherto escaped notice, are not without interest.''} Another
early work \cite{43} already contains a particular example of a
(self)distributive quasigroup:
$$\vbox{\offinterlineskip
\halign{\strut\hskip 4 pt # &\vrule height 11 pt \quad #\quad &
#\quad & #\,\cr
\phantom{0} & 0 & 1 & 2\cr\noalign{\hrule}
0 & 0 & 2 & 1\cr
1 & 2 & 1 & 0\cr
2 & 1 & 0 & 2\cr}}$$
This quasigroup is necessarily non-associative and plays a
principal r\^ole in the structure theory of distributive (or,
more generally, trimedial) quasigroups (see e.g.\ \cite{2},
\cite{3}, \cite{5}, \cite{6}, \cite{19}, \cite{35} and
\cite{48}).

The first article fully devoted to selfdistributivity
is (perhaps) \cite{11} (see also \cite{49} and \cite{32}) where,
among others, normal subquasigroups are studied and an attempt is
made to show that every minimal subquasigroup of a (finite)
distributive quasigroup is normal (see also \cite{15}). Actually,
the latter assertion is not true. All non-medial symmetric
distributive quasigroups (alias non-desarguesian planarily affine
triple systems) serve as counterexamples and first constructions
of these can be found in \cite{9} and \cite{17}. However, the
paper \cite{11} may be regarded as the starting point for the
investigation of normality problems in distributive quasigroups.

Hamiltonian groups (i.e., (non-commutative) groups having only
normal subgroups) were described (and named after W.\ R.\
Hamilton) by R.\ Dedekind in \cite{13} and it was shown in
\cite{38} that a similar description takes place for hamiltonian
Moufang loops, too. Furthermore, all subquasigroups of medial
quasigroups (i.e., quasigroups satisfying the identity
$(xy)(uv)\bumpeq(xu)(yv)$) are normal. That is, these quasigroups
are hamiltonian. (Notice that abelian groups are included in
hamiltonian structures in this paper -- not usual, but technically
advantageous.)

A thorough treatment (remarkable also for epic width) on cancellative
distributive groupoids was written by J.-P.\ Soublin (\cite{48}).
Section IV.9 of \cite{48} is devoted to normal subquasigroups of
distributive quasigroups and, among others, it is shown that
every hamiltonian symmetric (i.e., satisfying the identities
$xy\bumpeq yx$ and $x(xy)\bumpeq y$) distributive quasigroup is medial.
Moreover, an open problem whether there exist non-medial
hamiltonian distributive quasigroups is formulated (\cite{48},
p.\ 175). The main aim of the present paper is to solve this
problem.

In \cite{42}, it is claimed that every hamiltonian quasigroup
which is either distributive or a $CH$-quasigroup (i.e., a
symmetric quasigroup satisfying the identity $(xx)(yz)\bumpeq
(xy)(xz)$), is medial. The proof is based on the
idea that if $H$ is a~subloop of a~commutative Moufang loop
$G$ and the subloop generated by $H$ and the centre of $G$ is
normal then $H$ is normal. However, this assertion is false,
any non-associative commutative Moufang loop nilpotent of class 2
serving as an easy counterexample (in this case, every subloop
containing the centre is normal and $G$ contains a non-normal
subloop). Moreover, 3.2 and 8.9 are examples of non-medial
hamiltonian distributive or $CH$-quasigroups, respectively.

A possible way how to construct non-medial hamiltonian
distributive quasigroups is suggested in \cite{22},
but the paper is almost
unreadable and much more has to be done. However, the basic idea
is working, and the problem is transferred, step by step, first
to commutative Moufang loops with operators, then to certain
ternary algebras and, finally, to some cocyclic modules.
Actually, the problem of finding non-medial hamiltonian
distributive quasigroups is equivalent to the construction of
(finite) cocyclic modules over the ring $\Bbb{Z}[x, x^{-1},
(1-x)^{-1}]$ that cannot be generated by less than three elements.
We recall that a cocyclic module is contained in the injective
hull of its simple essential socle, so a good understanding of
the injective hull of simple modules and its submodules is
necessary to solve the problem.

After \cite{37}, \cite{30} and \cite{16}, if $R$
is a commutative noetherian ring then
the structure of some of the indecomposable injective modules over
$R[x]$, and hence over a~localization of $R[x]$, can be described in
terms of modules of divided powers over the  indecomposable
injectives of $R$. This is the case for  the injective hull of the
simple modules over the ring $R[x]=\Bbb{Z}[x]$. Since the
indecomposable injective modules over $\Bbb{Z}$ are also well
known, a detailed study of the modules  of divided powers  and
some of their finite submodules gives us the desired examples of
cocyclic modules, which, after the proper translation, allows us to
construct our examples of non-medial hamiltonian distributive
quasigroups in a completely explicit way.

Modules of divided powers, also called Macaulay modules, were
first known in \cite{33} and, as mentioned before, they
are important in connection to the description of injective
modules. If $K$ is any field and $K(x)$ is the field of fractions
of $K[x]$ then $M=K(x)/K[x]_{(x)}$ is an indecomposable injective
module with simple essential socle $K$. Note that
$M$  has $\{x^{-n} + K[x]_{(x)}\}_{n\in \Bbb{N}}$ as $K$-basis.
The modules of divided powers can be seen as an abstraction of
the structure of $M$ to the general setting of modules over a
polynomial ring.

The following text is divided into Sections 1 -- 12. Basic
notions are introduced in Section 1. Section 2 is devoted to
normal subquasigroups and Section 3 contains two examples, the
second one being (in view of 12.8) the solution of our problem
(it could be interesting to show the required properties
of this example directly, probably using a~computer). In Sections
4 and 5, some basic properties of commutative Moufang loops and
quasimodules (i.e., commutative Moufang loops with operators)
are investigated. Section 6 deals with ternary representations of
quasimodules. Section 7 is devoted to the connection between
hamiltonian quasimodules and certain cocyclic modules. In
Sections 8 and 9, (hamiltonian) trimedial and distributive
quasigroups, respectively, are studied and a transfer to
quasimodules is presented. Sections 10 and 11 are devoted to
modules of divided powers. In Section 12, a synthesis of the
preceding parts is made. The initial
problem is solved, but a~complete description of non-medial
hamiltonian distributive quasigroups is far from being finished.

\head 1. Preliminaries\endhead

\subhead\nofrills 1.1 (Quasigroups)\,\,\endsubhead
A non-empty set $Q$ equipped with a binary operation is said to be
a {\sl quasigroup} if for all $a,b\in Q$ there exist uniquely
determined $x,y\in Q$ such that $ax=b=ya$.
A quasigroup with a neutral element (a unit) is a {\sl loop}.

A quasigroup $Q$ is called
\roster
\item"--" {\sl medial} if $(ax)(yb)=(ay)(xb)$ for all $a,b,x,y\in
Q$;
\item"--" {\sl trimedial} if every subquasigroup of $Q$ generated
by at most three elements is medial;
\item"--" {\sl left (right) distributive} if $x(ab)=(xa)(xb)$ (\,
$(ab)x=(ax)(bx)$\,) for all $a,b,x\in Q$;
\item"--" {\sl distributive} if $Q$ is both left and right
distributive;
\item"--" {\sl symmetric} if $ax=xa$ and $x(xa)=a$ for all $a,x\in
Q$;
\item"--" a {\sl $CH$--quasigroup} if $Q$ is symmetric and
$(xx)(ab)=(xa)(xb)$ for all\linebreak $a,b,x\in Q$.
\endroster

Every distributive quasigroup is idempotent and trimedial
(\cite{2}). Every $CH$--quasigroup is trimedial (\cite{35}). A
reader is referred to \cite{2}, \cite{3}, \cite{11}, \cite{12},
\cite{15}, \cite{18}, \cite{35}, \cite{40}, \cite{48}, \cite{49}
for many useful prerequisites concerning (distributive, medial,
etc.) quasigroups.

\subhead\nofrills 1.2 (Commutative Moufang loops)\,\,\endsubhead
Let $Q$ be a loop satisfying the equation $(xx)(ab)=(xa)(xb)$.
Substituting $a=1_Q$, we get $(xx)b=x(xb)$ and, setting $b=1_Q$,
we get $(xx)a=(xa)x$. Now, if $a=b$, then $x(xa)=(xa)x$ and it
follows easily that $Q$ is commutative. Such a loop $Q$ is called
a {\sl commutative Moufang loop}. All the details concerning
commutative Moufang loops needed in the sequel may be found in
\cite{10}.

\subhead\nofrills 1.3 (Rings and modules)\,\,\endsubhead
In what follows, $\R$ stands for a (non--trivial) commutative and
associative noetherian ring with unit and modules are unitary
  $\R$--modules with scalars written on the left. Furthermore, we assume that there exists a
(ring) homomorphism $\F$ of $\R$ onto the three--element field
$\Bbb Z_3=\{0,1,2\}$ of integers modulo 3 and we put $\I=\Ker(\F)$.
Clearly, $\I$ is a maximal ideal of $\R$ and the simple
(three-element) factormodule $_\R\R/\I$ will be denoted by $\P$.
As concerns
various further pieces of information on general rings and
modules, a reader is referred to \cite{1}, \cite{8} and \cite{50}
and to \cite{36} for more specific information on the commutative
noetherian setting. A very nice reference for injective modules is
the book \cite{44}. The injective modules we study are, in fact,
artinian; for some of the results on artinian modules over
commutative ring we need the reference \cite{14}.

\subhead\nofrills 1.4 (Quasimodules)\,\,\endsubhead
By a {\sl quasimodule} we mean a commutative Moufang loop $Q(+)$
(usually denoted additively with neutral element 0) together with a scalar
multiplication $\R\times Q\to~Q$ such that the usual unitary
$\R$--modules equations are satisfied (i.e., $r(x+y)=rx+ry$,
$(r+s)x=rx+sx$, $(rs)x=r(sx)$, $1x=x$ and $0x=0$ for all
$r,s\in\R$ and $x,y\in Q$) and, moreover,
$rx+(y+z)=(rx+y)+z$ for all $r\in\I$ and $x,y,z\in Q$. The
quasimodule $Q$ is said to be {\sl primitive} if $\I Q=0$.
Obviously, if $Q$ is primitive then every subloop of $Q(+)$ is a
subquasimodule. See
\cite{21}, \cite{22}, \cite{25}, \cite{26}, \cite{28} and
\cite{29} for more on quasimodules.

\subhead\nofrills 1.5 (Ternary algebras)\,\,\endsubhead
By a {\sl ternary algebra} we mean a module $_\R A$ together with a
trilinear mapping $\tau:A^{(3)}\to A$ such that the following
equations are satisfied:
$$
\gather
\I\tau=0\,;\tag T0\\
\tau(x,x,y)=0\,;\tag T1\\
\tau(\tau(x,y,z),u,v)=0\,;\tag T2\\
\tau(u,v,\tau(x,y,z))=0\,.\tag T3
\endgather
$$

If $A=A(+,rx,\tau)$ is a ternary algebra then we put
$$
\overline{\tau}(x,y,z)=\tau(x,y,z)+\tau(y,z,x)+\tau(z,x,y)
$$
for all $x,y,z\in A$. Further,
$$
\An(a)=\{a\in A\,|\,\tau(a,x,y)=\tau(x,y,a)=0\text{ for all
}x,y\in A\}\,.
$$

\head 2. Normal subquasigroups\endhead

An equivalence $r$ defined on a quasigroup $Q$ is said to be a
{\sl normal congruence} of $Q$ if the following three conditions
are satisfied for all $a,b,c,d\in Q$:
$$
\gather
(a,b)\in r\text{ and }(c,d)\in r\Rightarrow (ac,bd)\in r\,;\tag
C1\\
(a,b)\in r\text{ and }(ac,bd)\in r\Rightarrow (c,d)\in r\,;\tag
C2\\
(c,d)\in r\text{ and }(ac,bd)\in r\Rightarrow (a,b)\in r\,.\tag
C3
\endgather
$$
(Note that both (C2) and (C3) follow from (C1) for a finite $Q$.)

\proclaim{2.1 Lemma} Let a subquasigroup $P$ of a quasigroup $Q$
be a block (or a class) of a normal congruence $r$ of $Q$. Then:
\roster
\widestnumber\item{}
\item"{\rm(i)}" Every block of $r$ is equal to a left coset $aP$ for
some $a\in Q$.
\item"{\rm(ii)}" Every block of $r$ is equal to a right coset $Pb$ for
some $b\in Q$.
\item"{\rm(iii)}" $(c,d)\in r\Leftrightarrow cP=dP\Leftrightarrow
Pc=Pd$.
\endroster
\endproclaim

\demo{Proof} Well known and easy.
\qed\enddemo

Now, a subquasigroup $P$ is said to be {\sl normal} in $Q$ if $P$
is a block of some normal congruence $r$ of $Q$; then, due to
2.1, $r$ is uniquely determined by $P$.

\proclaim{2.2 Lemma} A subquasigroup $P$ of a left distributive
quasigroup $Q$ is normal if and only if the following condition
is satisfied:
\roster
\widestnumber\item{}
\item"{\rm(C4)}" If $a,b,x,y,z\in Q$ are such that
$(xa)(yb)=z(ab)$ and if any two of the elements $x,y,z$ are in
$P$ then the remaining one is in $P$.
\endroster
\endproclaim

\demo{Proof} Assume first that $P$ is a block of a normal
congruence $r$ of $Q$. If $x,y\in P$ then $(xb,yb)\in r$,
$((xa)(xb),(xa)(yb))\in r$ and, since $(xa)(xb)=x(ab)$, we get
$(x(ab),z(ab))\in r$. Now, $(x,z)\in r$ by (C3) and consequently
$z\in P$. The other cases are similar.

Now, assume that (C4) is true and define a binary relation $r$ on
$Q$ by $(a,b)\in r\Leftrightarrow Pa=Pb$. If $(a,b)\in r$,
$(c,d)\in r$ and $x\in P$ then $x(ac)=(xa)(xc)=(yb)(zd)=w(bd)$
for suitable $x,y,z\in P$ and $w\in Q$. Using (C4), we get $w\in
P$ and the inclusion $Pac\subseteq Pbd$ follows. Quite similarly,
$Pbd\subseteq Pac$ and hence $(ac,bd)\in r$ and (C1) is verified.
The conditions (C2) and (C3) may be checked in a similar way.
Thus $r$ is a normal congruence and $P$ is among the blocks of
$r$ due to the definition of $r$ and the fact that $P$ is a
subquasigroup of $Q$.
\qed\enddemo

\proclaim{2.3 Lemma} Let $P$ be a subquasigroup of a left
distributive quasigroup $Q$. Then:
\roster
\widestnumber\item{}
\item"{\rm(i)}" $P\cdot ab\subseteq Pa\cdot Pb$ and
$|P|\le|Pa\cdot Pb|\le|P|^2$ for all $a,b\in Q$.
\item"{\rm(ii)}" If $P$ is finite then $P$ is normal in $Q$ if
and only if $|P|=|Pa\cdot Pb|$ for all $a,b\in Q$.
\endroster
\endproclaim

\demo{Proof} (i) For every $x\in P$, $x\cdot ab=xa\cdot xb\in
Pa\cdot Pb$.
\newline
(ii) Combine (i) and 2.2.
\qed\enddemo

A quasigroup $Q$ is called {\sl simple} if $Q$ is non--trivial
and $\id_q$, $Q\times Q$ are the only normal congruences of $Q$.

A quasigroup $Q$ is called {\sl hamiltonian} if every
subquasigroup is normal in $Q$. Clearly, the class of hamiltonian
quasigroups is closed under taking subquasigroups and
factorquasigroups. Hamiltonian groups serve as first examples of
hamiltonian quasigroups and the next basic result is almost
immediate (as in the subcase of abelian groups).

\proclaim{2.4 Proposition} Every medial quasigroup is
hamiltonian.
\endproclaim

\demo{Proof} If $P$ is a subquasigroup of a medial quasigroup $Q$
then we define a relation $r$ on $Q$ by $(a,b)\in
r\Leftrightarrow Pa=Pb$. Using the medial law, it is
straightforward and easy to show that $r$ is a normal congruence
of $Q$ and $P$ is one of the blocks.
\qed\enddemo

\flushpar{\bf 2.5} {\smc Remark.} (i) Because of technical
reasons, we prefer to include abelian groups into the class of
hamiltonian groups.\newline
(ii) Hamiltonian loops were studied in \cite{38} (see also
\cite{10}).\newline
(iii) Quasigroups linear over abelian groups (see \cite{23},
\cite{24}) are hamiltonian and play the r\^ole of
abelian groups among hamiltonian quasigroups.

\proclaim{2.6 Proposition} Let $Q$ be a left distributive
quasigroup and let $a\in Q$ be any element. The following
conditions are equivalent:
\roster
\item"{\rm(i)}" For every $x\in Q$, $x\ne a$, the subquasigroup
generated by the elements $a,x$ is normal in $Q$.
\item"{\rm(ii)}" Every two--generated subquasigroup is normal in
$Q$.
\item"{\rm(iii)}" $Q$ is hamiltonian.
\endroster
\endproclaim

\demo{Proof} (i) $\Rightarrow$ (ii). Let $u,v\in Q$, $u\ne
v$, and $P=\langle u,v\rangle$. Since $Q$ is a quasigroup, there
exist $b,x\in Q$ such that $ba=u$ and $bx=v$. Then $a\ne x$ and
$S=\langle a,x\rangle$ is a normal subquasigroup of $Q$ by (i).
On the other hand, the left translation $L_b:y\mapsto by$ is an
automorphism of $Q$ and hence $P=L_b(S)$ is normal in $Q$,
too.
\newline
(ii) $\Rightarrow$ (iii). We are going to check the condition (C4)
for a subquasigroup $P$ of $Q$ (see 2.2). Let $a,b,x,y,z\in Q$ be
such that $(xa)(yb)=z(ab)$ and $x,y\in P$ (the other two cases
are similar). If $P_1=\langle x,y\rangle$ then $P_1\subseteq P$
and $P_1$ is either trivial or two--generated. Thus $P_1$ is
normal in $Q$, $z\in P_1$ by (C4) for $P_1$ and, finally, $z\in
P$.
\newline
(iii) $\Rightarrow$ (i). This implication is trivial.
\qed\enddemo

\proclaim{2.7 Corollary} Let $Q$ be a finite distributive
quasigroup and let $a\in Q$. Then $Q$ is hamiltonian if and only
if $|\langle x,a\rangle|=|\langle x,a\rangle a_1\cdot\langle
x,a\rangle b_1|$ for all $x,a_1,b_1\in Q$, $x\ne a$, $a_1\ne
b_1$.
\qed\endproclaim

\head 3. Two examples\endhead

\flushpar{\bf 3.1} (\cite{9}, \cite{17})
Put $\Cal D_1=\Bbb Z_3\times\Bbb Z_3\times\Bbb Z_3\times\Bbb Z_3$
and define an operation $\bigtriangleup$ on $\Cal D_1$ by
$$\split
a\bigtriangleup b=(&2a(0)+2b(0)+a(1)a(3)b(2)+2a(2)a(3)b(1)+
2a(1)b(2)b(3)\\&+a(2)b(1)b(3),2a(1)+2b(1),2a(2)+2b(2),2a(3)+2b(3))
\endsplit$$
for all $a=(a(0),a(1),a(2),a(3))\in\Cal D_1$ and
$b=(b(0),b(1),b(2),b(3))\in\Cal D_1$. One may check easily that $\Cal
D_1(\bigtriangleup)$ is a symmetric distributive quasigroup (and hence
a $CH$--quasigroup) of order 81. On the other hand, if
$x=(0,0,0,0)$, $y=(0,1,0,0)$, $u=(0,0,1,0)$ and $v=(0,0,0,1)$
then
$$
(x\bigtriangleup y)\bigtriangleup(u\bigtriangleup v)=(1,1,1,1)\ne
(2,1,1,1)=(x\bigtriangleup u)\bigtriangleup(y\bigtriangleup v)\,,
$$
and so $\Cal D_1(\bigtriangleup)$ is not medial.

Furthermore,
$\Cal P=\{(0,0,0,0),(0,0,0,1),(0,0,0,2)\}$ is a three--element
subquasi\-group of $\Cal D_1(\bigtriangleup)$ and the set $(\Cal
P\bigtriangleup(0,0,1,0))\bigtriangleup(\Cal
P\bigtriangleup(0,1,0,0))$ contains just 9 elements. In view of
2.7, $\Cal P$ is not normal in $\Cal D_1(\bigtriangleup)$.

\subhead\nofrills 3.2\,\,\endsubhead
Put $\Cal D_2=\Bbb Z_{27}\times\Bbb Z_9\times\Bbb Z_3$ and define
an operation $\bigtriangledown$ on $\Cal D_2$ by
$$
\split
&a\bigtriangledown
b=(26a(0)+3a(1)+2b(0)+24b(1)+18a(0)a(2)b(1)+9a(0)b(1)b(2)\\
&+18a(1)b(0)b(2)+9a(1)a(2)b(0),
8a(1)+3a(2)+2b(1)+6b(2),2a(2)+2b(2))
\endsplit
$$
for all $a=(a(0),a(1),a(2))\in\Cal D_2$ and
$b=(b(0),b(1),b(2))\in \Cal D_2$. Again, a tedious but
straightforward calculation shows that $\Cal D_2$ is a
distributive quasigroup of order 729 and $\Cal D_2$ is not
medial, since
$$(x\bigtriangledown y)\bigtriangledown(u\bigtriangledown v)=(7,5,1)\ne
(25,5,1)=(x\bigtriangledown u)\bigtriangledown(y\bigtriangledown v)\,,$$
where $x=(0,0,0)$, $y=(1,0,0)$, $u=(0,1,0)$, $v=(0,0,1)$. Finally,
using 2.6 with $a=(0,0,0)$ and 2.2 or 2.7, one may also (at least in
principle) show that $\Cal D_2(\bigtriangledown)$ is hamiltonian.
Nevertheless, this property is a consequence of 12.8.

\head 4. Commutative Moufang loops\endhead

Let $Q=Q(+)$ be a commutative Moufang loop, the operation
being denoted additively. The set
$$
\Z(Q)=\{a\in Q\,\,|\,\,(a+x)+y=a+(x+y)\text{ for all }x,y\in Q\}
$$
is a normal subloop of $Q$, called the {\sl centre} of $Q$. The
loop $Q$ is said to be {\sl nilpotent of class at most 0} if it
is trivial, {\sl of class at most 1} if it is an (abelian) group and
{\sl of class at most} $n\ge 2$ if the factorloop $Q/\Z(Q)$ is nilpotent
of class at most $n-1$. Further, $Q$ is {\sl nilpotent of class} $n$ if
it nilpotent of class at most $n$ and is not nilpotent of class
at most $n-1$. The smallest normal subloop $P$ of $Q$
such that the corresponding factorloop $Q/P$ is associative is
called the {\sl associator subloop} of $Q$ and is denoted by
$P=\A(Q)$ in the sequel. For all $a,b,c\in Q$, the element
$[a,b,c]=((a+b)+c)-(a+(b+c))$ is called the {\sl associator} of
$a,b,c$. Clearly, $\A(Q)$ is just the subloop generated by all
associators.

\proclaim{4.1 Proposition} {\rm(\cite{10})}
\roster
\runinitem"{\rm(i)}" Both $\A(Q)$ and
$Q/\Z(Q)$ are 3-elementary loops {\rm (i.e., they satisfy the
equation $3x=0$)}.
\widestnumber\item{}
\item"{\rm(ii)}" $Q$ is diassociative {\rm (i.e., any two elements
generate a subgroup)}.
\item"{\rm(iii)}" If $a,b,c\in Q$ are such that $a+(b+c)=(a+b)+c$ then
these elements generate a subgroup.
\item"{\rm(iv)}" If $Q$ is generated by $n\ge 2$ elements then it is
nilpotent of class at most $n-1$.
\item"{\rm(v)}" $Q$ is locally nilpotent.
\item"{\rm(vi)}" If $Q$ is simple then it is an abelian group of finite
prime order.
\item"{\rm(vii)}" If $Q$ is finite and not associative then
the order of $Q$ is divisible by 81.\qed
\endroster
\endproclaim

\flushpar{\bf 4.2} {\smc Remark.} It is proven in \cite{4},
\cite{34} and \cite{45} that the free commutative Moufang loop of
rank $n\ge 2$ is nilpotent of class $n-1$.\par\medskip

A transformation $f$ of $Q$ is said to be {\sl central} (more
precisely, $n$--{\sl central}) if there exists $n\in\Bbb Z$ such
that $f(x)+nx\in\Z(Q)$ for every $x\in Q$.

\proclaim{4.3 Lemma} {\rm(\cite{29})} Let $f$ be a transformation
of $Q$. Then:
\roster
\widestnumber\item{}
\item"{\rm(i)}" If $m,k\in\Bbb Z$, $0\le m\le 2$ and $n=3k+m$ then
$f$ is $n$--central if and only if it is $m$--central.
\item"{\rm(ii)}" If $Q$ is non--associative and $f$ is central
then there is just one $r\in\Bbb Z$ such that $0\le r\le 2$ and
$f$ is $r$--central.\qed
\endroster
\endproclaim

\proclaim{4.4 Lemma} {\rm(\cite{29})} Let $f$ and $g$ be
endomorphisms of $Q$ such that $f$ and $g$ are $m$--central and
$n$--central, respectively. Then:
\roster
\widestnumber\item{}
\item"{\rm(i)}" $fg$ is $(-mn)$--central.
\item"{\rm(ii)}" $f+g$ is an $(m+n)$--central endomorphism.
\item"{\rm(iii)}" If $f$ is an automorphism then $f^{-1}$ is
$m$--central.\qed
\endroster
\endproclaim

\flushpar{\bf 4.5} {\smc Remark.} Assume that $Q$ is not
associative. By 4.3 and 4.4, the set $\Cend(Q)$ of central
endomorphisms of $Q$ is an associative ring with unit and, for
every $f\in\Cend(Q)$, there is a uniquely determined
$\varPhi(f)\in\Bbb Z_3$ such that $f$ is
$(-\varPhi(f))$--central. Now, the mapping
$\varPhi:\Cend(Q)\to\Bbb Z_3$ is a projective ring homomorphism and
$\Ker(\varPhi)=\{f\,|\,f(Q)\subseteq\Z(Q)\}$.

\head 5. Quasimodules\endhead

Throughout this section, let $Q$ be a quasimodule, the underlying
commutative Moufang loop of $Q$ being denoted by $Q(+)$. A
subquasimodule $P$ of $Q$ is {\sl normal} if $P$ is a block of a
congruence of $Q$. If $Q$ is finitely generated then $\gen(Q)$ is
the smallest number of generators of $Q$.

\proclaim{5.1 Proposition} {\rm(\cite{21})}
\roster
\runinitem"{\rm(i)}" A subquasimodule $P$ of $Q$ is normal if and only
if $P(+)$ is a normal subloop of $Q(+)$.
\widestnumber\item{}
\item"{\rm(ii)}" $\A(Q)$ is a normal primitive subquasimodule of
$Q$ and $Q/\A(Q)$ is a module.
\item"{\rm(iii)}" $\Z(Q)$ is a normal submodule of $Q$ and
$Q/\Z(Q)$ is a primitive quasimodule.
\item"{\rm(iv)}" If $P$ is a subquasimodule of $Q$ such that
either $\A(Q)\subseteq P$ or $P\subseteq\Z(Q)$ then $P$ is a
normal subquasimodule of $Q$.
\item"{\rm(v)}" For all $x,y\in Q$, the set $\R x+\R y$ is a
submodule of $Q$ and it is just the subquasimodule generated by
$x,y$.
\item"{\rm(vi)}" If $a+(b+c)=(a+b)+c$ for some $a,b,c\in Q$ then
the subquasimodule generated by these elements is a submodule.
\item"{\rm(vii)}" If $Q$ is simple (i.e., if $Q\ne0$ and
$0,Q$ are the only normal subquasimodules of $Q$) then $Q$ is a
module (and hence $0,Q$ are the only submodules of $Q$).\qed
\endroster
\endproclaim

A {\sl preradical} $\varrho$ (for the category of quasimodules)
is a subfunctor of the identity functor, i.e., $\varrho$ assigns
to each quasimodule $Q$ its subquasimodule $\varrho(Q)$ in such a
way that $f(\varrho(Q))\subseteq \varrho(P)$ whenever $P$ is a
quasimodule and $f:Q\to P$ is a homomorphism. Obviously,
$\varrho(Q)$ is a normal subquasimodule of $Q$. A preradical
$\varrho$ is said to be {\sl hereditary} if
$\varrho(P)=P\cap\varrho(Q)$ for every quasimodule $Q$ and its
subquasimodule $P$, and it is said to be a {\sl radical} if
$\varrho(Q/\varrho(Q))=0$ for every quasimodule $Q$.

Let $\varrho$ be a preradical. For every quasimodule $Q$ and
every ordinal $\alpha$ we put $^0\!\varrho(Q)=0$,
$^{\alpha+1}\!\varrho(Q)/^{\alpha}\!\varrho(Q)=
\varrho(Q/^{\alpha}\!\varrho(Q))$,
$^{\alpha}\!\varrho(Q)=\bigcup_{\beta<\alpha}\,^{\beta}\!\varrho(Q)$
for $\alpha$ limit, and
$\widehat{\varrho}(Q)=\bigcup\,^{\alpha}\!\varrho(Q)$. Then
$^{\alpha}\!\varrho$, $\widehat{\varrho}$ are preradicals which are
hereditary if $\varrho$ is hereditary, and $\widehat{\varrho}$ is
the least radical containing $\varrho$ (see \cite{21}).

Let $\K(Q)$ denote the greatest primitive subquasimodule of $Q$
and $\S(Q)=\roman{Soc}(Q)$ the socle of $Q$ (i.e., the
subquasimodule generated by all minimal submodules). Then we have
$\A(Q)\subseteq\K(Q)\subseteq\S(Q)$. Moreover, both $\K$ and $\S$
are hereditary preradicals for the category of quasimodules. Now,
$\widehat{\K}$ and $\widehat{\S}$ will denote the smallest
hereditary radical containing $\K$ and $\S$, respectively.

\proclaim{5.2 Lemma} Let $P$ be a subquasimodule of $Q$. Then:
\roster
\widestnumber\item{}
\item"{\rm(i)}" If $P\cap\Z(Q)=0$ then $P\subseteq\K(Q)$.
\item"{\rm(ii)}" If $P$ is cyclic and $P\cap\Z(Q)=0$ then either
$P=0$ or $P\simeq\P$.
\item"{\rm(iii)}" If $P\ne 0$ is normal and cyclic then
$P\cap\Z(Q)\ne0$ and, moreover, if $P$ is simple then $P\subseteq\Z(Q)$.
\endroster
\endproclaim

\demo{Proof} (i) $P$ is isomorphic to a subquasimodule of
$Q/\Z(Q)$, and hence $P$ is primitive by 5.1(iii).\newline
(ii) By (i), $P$ is a cyclic $\K$--torsion module.\newline
(iii) Assume on the contrary that $P\cap\Z(Q)=0$. By (ii), $P$
contains just three elements, so that $P=\{0,a,-a\}$, $a\ne 0$.
Now, for $x,y\in Q$, put $z=[x,y,a]=((x+y)+a)-(x+(y+a))$. Since
$P$ is normal in $Q$, we have $z\in P$. If $z=a$ then
$(x+y)+a=(x+(y+a))+a$, $x+y=x+(y+a)$, $y=y+a$ and $a=0$, a
contradiction. If $z=-a$ then
$x+(y+a)=((x+y)+a)+a=(x+y)+2a=(x+a)+(y+a)$, $x=x+a$ and $a=0$,
again a contradiction. Thus $z=0$ and $(x+y)+a=x+(y+a)$. This
means $a\in\Z(Q)$, a final contradiction.
\qed\enddemo

\proclaim{5.3 Lemma} {\rm(\cite{25})} Assume that $Q$ is not
associative and $\gen(Q)=3$. Then:
\roster
\widestnumber\item{}
\item"{\rm(i)}" $\A(Q)\simeq\P$ and $Q/\Z(Q)\simeq\P^{(3)}$.
\item"{\rm(ii)}" If $P$ is a proper subquasimodule of $Q$ with
$\Z(Q)\subseteq P$ then $P$ is a module.
\item"{\rm(iii)}" If $P$ is a non--associative subquasimodule of
$Q$ then $\A(Q)\subseteq P$ and $P$ is normal in $Q$.
\item"{\rm(iv)}" If $Q$ is $\widehat{\K}$--torsion then every
proper subquasimodule of $Q$ is a module.\qed
\endroster
\endproclaim

The quasimodule $Q$ is said to be {\sl nilpotent of class at most}
$n\ge 0$ if so is the underlying commutative Moufang loop $Q(+)$.

\proclaim{5.4 Proposition} {\rm(\cite{21})} Assume that $Q$ is
finitely generated. Then:
\roster
\widestnumber\item{}
\item"{\rm(i)}" If $\gen(Q)=n\ge 2$ then $Q$ is nilpotent of class
at most $n-1$.
\item"{\rm(ii)}" $Q$ is noetherian {\rm(i.e., every
subquasimodule of $Q$ is finitely generated)}.\qed
\endroster
\endproclaim

\proclaim{5.5 Lemma} Assume that $Q$ is subdirectly irreducible
and nilpotent of class 2. Then $\A(Q)\simeq\P$ and every proper
factorquasimodule of $Q$ is a module.
\endproclaim

\demo{Proof} Since $0\ne\A(Q)\subseteq\Z(Q)$, $\A(Q)$ is a
subdirectly irreducible primitive module and consequently
$\A(Q)\simeq\P$.
\qed\enddemo

\proclaim{5.6 Proposition} Assume that $Q$ is
$\widehat{\K}$--torsion and not associative. The following
conditions are equivalent:
\roster
\item"{\rm(i)}" $Q$ is subdirectly irreducible and $\gen(Q)=3$.
\item"{\rm(ii)}" Every proper factorquasimodule as well as every
proper subquasimodule of $Q$ is a module.
\endroster
\endproclaim

\demo{Proof} (i) $\Rightarrow$ (ii). By 5.4(i), $Q$ is nilpotent of
class 2 and we can use 5.5 and 5.3(iv).\newline
(ii) $\Rightarrow$ (i). Since $Q$ is not associative, we have
$(a+b)+c\ne a+(b+c)$ for some $a,b,c\in Q$ and it is clear that
$Q$ is generated by these elements. Thus $\gen(Q)=3$. The fact
that $Q$ is subdirectly irreducible is also clear.
\qed\enddemo

\proclaim{5.7 Lemma} Assume that $Q$ is finitely generated and
let $P$ be a (proper) maximal subquasimodule of $Q$. Then $P$ is
normal in $Q$.
\endproclaim

\demo{Proof} We shall proceed by induction on the nilpotence
class $n$ of $Q$ (see 5.4(i)). First, the result is clear for
$n\le 1$ and if $\Z(Q)\subseteq P$ then $P/\Z(Q)$ is normal in
$Q/\Z(Q)$ by induction. Thus $P$ is normal in $Q$ in this case and
we may assume that $\Z(Q)\nsubseteq P$. But then $Q=P+\Z(Q)$ and
it is easy to check directly that $P$ is normal in $Q$.
\qed\enddemo

If $Q$ is finitely generated then $\J(Q)$ will denote the
intersection of all maximal submodules of $Q$.

\proclaim{5.8 Proposition} Assume that $Q$ is finitely generated.
Then:
\roster
\widestnumber\item{}
\item"{\rm(i)}" $\J(Q)$ is a normal subquasimodule of $Q$ and
$\A(Q)\subseteq\J(Q)$.
\item"{\rm(ii)}" $\gen(Q)=\gen(Q/\J(Q))=\gen(Q/\A(Q))$.
\item"{\rm(iii)}" If $Q$ is $\widehat{\K}$--torsion then
$Q/\J(Q)$ is primitive and $|Q/\J(Q)|=3^{\gen(Q)}$.
\endroster
\endproclaim

\demo{Proof} By 5.7 and 5.1(vii), $\J(Q)$ is a normal
subquasimodule and $\A(Q)\subseteq\J(Q)$. The inequalities
$n=\gen(Q)\ge\gen(Q/\A(Q))\ge\gen(Q/\J(Q))$ are clear. Now, let
$N$ be a generator set of $Q$, $|N|=n$, and let $M$ be a subset
of $Q$ such that $Q/\J(Q)$ is generated by $M/\J(Q)$. We claim
that $Q$ is generated by $M$.

Assume the contrary and consider a subset $N_1$ of $N$ maximal
with respect to the property that $Q$ is not generated by $M\cup
N_1$. Then $N_1\ne N$ and we take $v\in N\setminus N_1$. Further,
consider a subquasimodule $V$ of $Q$ maximal with respect to
$M\cup N_1\subseteq V$ and $v\notin V$. It is easy to see that $V$
is a maximal subquasimodule of $Q$, and hence $\J(Q)\subseteq V$
and $V=Q$, a contradiction.

We have shown that $M$ generates $Q$ and it follows easily that
$\gen(Q/\J(Q))=\gen(Q)$.

Now, finally, assume that $Q$ is $\widehat{\K}$--torsion. Then
every simple factor of $Q$ is a~copy of $\P$ and $Q/\J(Q)$ is a
primitive module which is a direct sum of $n$ copies of $\P$.
\qed\enddemo

\proclaim{5.9 Lemma} If $Q$ is finitely generated and
$\widehat{\K}$--torsion then $\I Q\subseteq\J(Q)\cap\Z(Q)$ and
$\gen(Q)=\gen(Q/\I Q)$.
\endproclaim

\demo{Proof} Use 5.8.
\qed\enddemo

\proclaim{5.10 Lemma} Assume that $Q$ is a primitive quasimodule
nilpotent of class at most 2 with $\gen(Q)=n$. Then $|Q|\le
3^{n+m}$, where $m=\binom{n}{3}$.\endproclaim

\demo{Proof} As $Q/\A(Q)$ is a primitive module, its additive
group is 3--elementary and every its subgroup is a submodule. By
5.8(ii), $\gen(Q/\A(Q))=n$, and consequently $|Q/\A(Q)|=3^n$. If
$n\le 2$ then $\A(Q)=0$. Now, let $n\ge 3$, $M=\{a_1,\dots,a_n\}$
be a generating set of $Q$, $N=\{\,[a_i,a_j,a_k]\,|\,1\le
i<j<k\le n\}$, and $P$ be the subquasimodule generated by $N$.
Then, $Q$ being nilpotent of class at most 2,
$P\subseteq\A(Q)\subseteq\Z(Q)$ and $P$ is a normal
subquasimodule of $Q$ by 5.1(iv). Denote by $f$ the natural
projection of $Q$ onto $Q/P$. Then $K=\{f(a_1),\dots,f(a_n)\}$
generates $Q/P$, $x+(y+z)=(x+y)+z$ for all $x,y,z\in K$ and $Q/P$
is associative by 5.1(vi). Thus $\A(Q)=P$, $\gen(\A(Q))\le|N|=m$
and $|\A(Q)|\le 3^m$, since $\A(Q)$ is a primitive module, too.
\qed\enddemo

\proclaim{5.11 Lemma} Let $P$ be a minimal submodule of $Q$,
$a,b\in Q$, $A=P+(a+b)$ and $B=(P+a)+(P+b)$. Then:
\roster
\widestnumber\item{}
\item"{\rm(i)}" $A\subseteq B$.
\item"{\rm(ii)}" If $P$ is normal in $Q$ then $A=B$.
\item"{\rm(iii)}" If $P$ is not isomorphic to $\P$ then $A=B$.
\item"{\rm(iv)}" If $A\ne B$ then $P\simeq\P$, $|A|=3$ and
$|B|=9$.
\endroster
\endproclaim

\demo{Proof} We may assume that $P$ is not normal in $Q$, $A\ne
B$ and $Q$ is generated by $P\cup\{a,b\}$. By 5.2(ii) and 5.4(i),
$P\simeq\P$ and $Q$ is nilpotent of class 2. Consequently, $|A|=3$
and $3\le|B|\le 9$. Now, let $(x+a)+(y+b)=(u+a)+(v+b)$ for some
$x,y,u,v\in P$, $x\ne u$, $y\ne v$. Then we have
$a+(r+b)=((x+a)+(y+b))-2x=((u+a)+(v+b))-2x=(a+s)+(b+t)$, where
$r=y-x$, $s=u-x$, $t=v-x$, $r\ne t$, $s\ne 0$ (the
subquasimodules $\langle P,a\rangle$ and $\langle P,b\rangle$ are
at most two--generated, and hence they are associative).
Furthermore, $\A(Q)\subseteq\Z(Q)$, and therefore
$(a+b)+r+\alpha=a+(r+b)=(a+s)+(b+t)=(a+b)+(s+t)+\beta$ for some
$\alpha,\beta\in\Z(Q)$. Then $\alpha+r=\beta+(s+t)$,
$\alpha-\beta\in P\cap\Z(Q)=0$ (since $P$ is minimal and not normal
in~$Q$), $\alpha=\beta$ and $r=s+t$. Thus $a+(b+r)=(a+s)+(b+t)$,
where $r=s+t$ and $s\ne 0$. If $r=0$ then
$(a+b)+2s=((a+s)+(b-s))+2s=(a+2s)+b$ and $Q=\langle a,b,2s\rangle$
is associative, a contradiction. Similarly if $t=0$. Finally, if
$r\ne0\ne t$ then $t=s$ and $a+(b+r)=(a+b)+2s=(a+b)+r$, again a
contradiction.
\qed\enddemo

\flushpar{\bf 5.12} {\smc Remark.} Let $P$ be a minimal submodule
of $Q$. Then, for every $a\in Q$, the subquasimodule $\langle
P,a\rangle$ is at most two--generated, and so it is a submodule
and $P+(P+a)=P+a$. By 5.11, $P$ is normal in $Q$ if and only if
$|P|\cdot|(P+a)+(P+b)|\ne 27$ for all $a,b\in Q$.

\head 6. Ternary representations of quasimodules nilpotent of
class at most 2\endhead

Throughout this section, the word {\sl quasimodule} always means
{\sl quasimodule nilpotent of class at most 2}.

\proclaim{6.1 Proposition} Let $A=A(+,rx,\tau)$ be a ternary
algebra. Then $\q(A)=A(\*,rx)$ is a quasimodule, where the
underlying commutative Moufang loop is defined by
$$
x\*y=x+y+\tau(x,y,x-y)
$$
for all $x,y\in A$. Moreover, $\Z(A(\*))=\{a\in
A\,|\,\overline{\tau}(a,x,y)=0\text{ for all }x,y\in A\}$ and
$\A(A(\*))$ is the subloop generated by $\Im(\overline{\tau})$.
\endproclaim

\demo{Proof} Clearly, 0 is the neutral element of $A(\*)$ and
$x\*y=x+y+\tau(x,y,x)+\tau(y,x,y)=y\*x$. Further,
$(x\*x)\*(y\*z)=2x+y+z+\tau(x,y,x)+\tau(x,y,y)+\tau(x,z,x)
+\tau(x,z,z)+\tau(y,z,y)+\tau(z,y,z)+\tau(x,y,z)+\tau(x,z,y)
=(x\*y)\*(x\*z)$ for all $x,y,z\in A$. If $x\*y=x\*z$ then
$y+\tau(x,y,x)+\tau(y,x,y)=v=z+\tau(x,z,x)+\tau(z,x,z)$, and
hence $\tau(x,y,x)=\tau(x,v,x)=\tau(x,z,x)$,
$\tau(y,x,y)=\tau(v,x,v)=\tau(x,z,x)$ and $y=z$. Finally, if
$z=y-x+\tau(x,y,x)+\tau(x,y,y)$ then $x\*z=y$. We have checked
that $A(\*)$ is a commutative Moufang loop. The opposite (or
inverse) element to $x$ is $-x$ and $x\ominus
y=x\*(-y)=x-y+\tau(y,x,x)+\tau(y,x,y)$. Now, for all $a,x,y\in
A$, we have $[a,x,y]=((a\*x)\*y)\ominus(a\*(x\*y))=
\overline{\tau}(a,x,y)$. Consequently,
$\Z(A(\*))=\{a\,|\,\overline{\tau}(a,x,y)=0\}$ and it is clear
that $\I A\subseteq\Z(A(\*))$.

For every $r\in\R$, $r^3-r\in\I$, and hence
$r(x\*y)=rx+ry+r\tau(x,y,x-y)=rx+ry+r^3\tau(x,y,x-y)=rx\*ry$.
Similarly, $(r+s)x=rx+sx=rx+sx+\tau(rx,sx,rx-sx)=rx\*sx$ and we
see that $A(\*,rx)$ is a quasimodule. It remains to show that
this quasimodule is nilpotent of class at most 2. However,
$((x\*y)\*z)\ominus(x\*(y\*z))\in\Z(A(\*))$ for all $x,y,z\in A$
and the rest is clear.
\qed\enddemo

Quasimodules $\q(A)$, $A$ being a ternary algebra, will be said to
have {\sl ternary representation}. Now, we are going to show
that every free quasimodule of finite rank has a ternary
representation.\par\medskip

\flushpar{\bf 6.2} Let $n\ge 2$, $m=\binom{n}{3}$, $q=n+m$, and
$F=F_n=\R^{(n)}\times\P^{(m)}$. Then $F$ is an $\R$--module
and the elements $a_1=(1,0,\dots,0),\dots,a_q=(0,\dots,0,1)$
form a~canonical set $M$ of generators of $F$. Let $K$ be
the set of ordered triples $(i,j,k)$, $1\le i<j<k\le n$, and let
$f:K\to\{1,\dots,m\}$ be a bijection. Now, define a~mapping
$\sigma:M^{(3)}\to F$ by $\sigma(a_i,a_j,a_k)=a_{n+f(\alpha)}$,
$\sigma(a_j,a_i,a_k)=2a_{n+f(\alpha)}$ for every
$\alpha=(i,j,k)\in K$ and $\sigma(a_i,a_j,a_k)=0$ for every
triple $(i,j,k)$ such that neither $(i,j,k)$ nor $(j,i,k)$ is in
$K$. Then this mapping $\sigma$ can be extended (in a~unique way)
to a trilinear mapping $\tau:F^{(3)}\to F$ and $F=F(+,rx,\tau)$
becomes a~ternary algebra. Now, consider the corresponding
quasimodule $\q(F)=F(\*,rx)$ (see 6.1).

\proclaim{6.2.1 Proposition} $\q(F)$ is a free quasimodule and the
set $N=\{a_1,\dots,a_n\}$ is a free basis of $\q(F)$. Moreover,
$\A(F(\*))=\P^{(m)}$, $\Z(F(\*))=\I^{(n)}\times\P^{(m)}$ and $\I
F=\I^{(n)}$.
\endproclaim

\demo{Proof} We have
$((a_i\*a_j)\*a_k)\ominus(a_i\*(a_j\*a_k))=a_{n+f(\alpha)}$ for
every $\alpha=(i,j,k)\in K$. Consequently, the quasimodule $q(F)$
is generated by $N$. The equalities $\A(F(\*))=\P^{(m)}$,
$\Z(F(\*))=\I^{(n)}\times\P^{(m)}$ and $\I F=\I^{(n)}$ are also
easy to check. It remains to show that $\q(F)$ is free over $N$.

Now, let $E=E(\*,rx)$ be the free quasimodule over $N$ and let
$\pi:E\to\q(F)$ be the (unique) projective quasimodule
homomorphism such that $\pi\restriction N=\id_N$. Then
$\pi(\A(E))=\A(\q(F))$, $\pi(\I E)=\I F$ and $\pi$ induces
projective homomorphisms $\varphi:E/\A(E)\to \q(F)/\A(\q(F))$ and
$\psi:E/\I E\to \q(F)/\I F$ such that $\varphi\lambda=\varrho\pi$ and
$\psi\mu=\nu\pi$, where $\lambda,\varrho,\mu,\nu$ are the
corresponding natural projections. Moreover,
$\varphi\lambda(N)=\varrho\pi(N)=\varrho(N)$ is a free basis of
the free module $\q(F)/\A(\q(F))$,
$|\lambda(N)|=|\varphi\lambda(N)|=|\varrho(N)|=|N|=n$ and we
conclude that $\varphi$ is an isomorphism and
$\Ker(\pi)\subseteq\Ker(\lambda)=\A(E)$. On the other hand,
$|\q(F)/\I F|=3^q$ and $E/\I E$ is a free primitive quasimodule of
rank $n$. By 5.10, $|E/\I E|\le 3^q$, and therefore
$\psi$ is also an isomorphism and $\Ker(\pi)\subseteq\Ker(\mu)=\I
E$. We have shown that $\Ker(\pi)\subseteq\A(E)\cap\I E$ and to
finish the proof it suffices to check that $\A(E)\cap\I E=0$.

First, take $a\in\I E$. Since $\I E\subseteq\Z(E)$ and $N$ is a
free basis of $E$, we have $a=r_1a_1\*\dots\*r_na_n$ for some
$r_1,\dots,r_n\in\R$. Now, if $a\in\A(E)$ then
$0=\lambda(a)=r_1\lambda(a_1)\*\dots\*r_n\lambda(a_n)$. But
$E/\A(E)$ is a free module over
$\{\lambda(a_1),\dots,\lambda(a_n)\}$, and hence
$r_1=\dots=r_n=0$ and $a=0$.
\qed\enddemo

\proclaim{6.2.2 Lemma} $\tau\restriction Z\times F\times F=\tau
\restriction F\times Z\times F = \tau\restriction F\times F\times
Z=0$, where $Z=\Z(F(\*))$.
\endproclaim

\demo{Proof} Obvious.
\qed\enddemo

\proclaim{6.2.3 Corollary} Every submodule of $\Z(F(\*))$ is an
ideal of the ternary algebra $F$.\qed
\endproclaim

\proclaim{6.3 Proposition} Every finite $\widehat{\K}$--torsion
quasimodule has ternary representation.
\endproclaim

\demo{Proof} Let $Q=Q(\*,rx)$ be a finite $\widehat{\K}$--torsion
quasimodule which is not a module. Then $n=\gen(Q)\ge 3$ and
there exists a projective homomorphism $\pi:\q(F)\to Q$ as it
follows from 6.2.1; put $G=\Ker(\pi)\*\I F$. Then $G$ is a normal
subquasimodule of $\q(F)$ and $\varphi(G)$ is a normal
subquasimodule of $H=\q(F)/\I F$, $\varphi:\q(F)\to H$ being the
natural projection. Moreover, $H_1=H/\varphi(G)\simeq Q/\I Q$
and, by 5.8 and 5.9, $n=\gen(Q/\I
Q)=\gen(H_1)=\gen(H_1/\A(H_1))$. Consequently, $H_1/\A(H_1)$ is a
primitive module of dimension $n$ and $|H_1/\A(H_1)|=3^n$. Since
$H_1/\A(H_1)\simeq H/L$, where $L=\varphi(G)\*\A(H)$, we have
also $|H/L|=3^n$. On the other hand, $|H/\A(H)|=3^n$, $\A(H)=L$,
$\varphi(G)\subseteq\A(H)$ and $\Ker(\pi)\subseteq\A(\q(F))\*\I
F=\Z(\q(F))$ (see 6.2.1). Now, by 6.2.3, $\Ker(\pi)$ is an~ideal
of the ternary algebra $F$ and it suffices to consider the
corresponding factoralgebra.
\qed\enddemo

\flushpar{\bf 6.4} {\smc Remark.} Using primary decompositions
(and filtered products), one may show that every finite
$\widehat{\S}$--torsion quasimodule has ternary representation
(and that every $\widehat{\S}$--torsion quasimodule is imbeddable
into a quasimodule with ternary representation).

\head 7. Hamiltonian quasimodules\endhead

A quasimodule $Q$ is said to be {\sl hamiltonian} if every
subquasimodule is normal in $Q$. Clearly, every module is
hamiltonian and the class of hamiltonian quasimodules is closed
under subquasimodules and factorquasimodules.

A quasimodule $Q$ is said to be {\sl cocyclic} if $\S(Q)$ is a
non--zero essential simple submodule of $Q$.

\proclaim{7.1 Proposition} Let $Q$ be a cocyclic quasimodule.
Then:
\roster
\widestnumber\item{}
\item"{\rm(i)}" $Q$ is subdirectly irreducible, hamiltonian and
nilpotent of class at most 2.
\item"{\rm(ii)}" If $Q$ is non--associative then
$\A(Q)=\S(Q)\simeq\P$ and every proper factorquasimodule of $Q$
is a module.
\endroster
\endproclaim

\demo{Proof} Clearly, $\S(Q)$ is the smallest non--zero normal
subquasimodule, and hence $Q$ is subdirectly irreducible. Now, we
may assume that $Q$ is not associative. Then
$\A(Q)\subseteq\S(Q)$ implies $\A(Q)=\S(Q)\simeq\P$ and
$\A(Q)\subseteq\Z(Q)$ (see 5.2(iii) and 5.5). Thus $Q$ is
hamiltonian and nilpotent of class at most 2.
\qed\enddemo

\proclaim{7.2 Lemma} If $Q$ is a non--associative hamiltonian
quasimodule then
$\A(Q)\subseteq\K(Q)\subseteq\S(Q)\subseteq\Z(Q)$.
\endproclaim

\demo{Proof} Use 5.2(iii).
\qed\enddemo

\proclaim{7.3 Corollary} Every hamiltonian primitive quasimodule
is a module.\qed
\endproclaim

\proclaim{7.4 Proposition} A non--zero quasimodule is cocylic if
and only if it is hamiltonian and subdirectly irreducible.
\endproclaim

\demo{Proof} Use 7.1.
\qed\enddemo

\proclaim{7.5 Proposition} Let $Q$ be a non--associative cocyclic
quasimodule. Then $Q$ is $\widehat{\K}$--torsion. Moreover, if
$Q$ is finitely generated then it is finite and $|Q|=3^n$ for
some $n\ge 4$.
\endproclaim

\demo{Proof} By 7.1, $\P\simeq\A(Q)=\S(Q)\subseteq\Z(Q)$. Since
$\R$ is a commutative noetherian ring, every hereditary radical
(for $\R$--Mod) is stable, and hence, in particular, $\Z(Q)$ is
$\widehat{\K}$--torsion. On the other hand, the factor $Q/\Z(Q)$
is primitive, thus being $\K$--torsion. Consequently, $Q$ is
$\widehat{\K}$--torsion.

Now, assume that $Q$ is finitely generated. By 5.4(ii), $Q$ is a
noetherian quasimodule, $Q$ has a finite $\K$--sequence and we
may restrict ourselves to the case when $Q$ is $\K$--torsion.
Then both $\Z(Q)$ and $Q/\Z(Q)$ are noetherian primitive modules,
thus being finite direct sums of copies of $\P$ and the rest is
clear.
\qed\enddemo

\proclaim{7.6 Proposition} Let $Q$ be a non--associative
hamiltonian quasimodule. Then there exist a subquasimodule $Q_1$
of $Q$ and a (normal) subquasimodule $Q_2$ of $Q_1$ such that the
factor $Q_3=Q_1/Q_2$ is non--associative, cocyclic, $\widehat{\K}$--torsion,
$\gen(Q_3)=3$ and $|Q_3|=3^n$ for some $n\ge 4$.
\endproclaim

\demo{Proof} Since $Q$ is not associative, there is a
non--associative subquasimodule $Q_1$ of $Q$ such that
$\gen(Q_1)=3$. Further, there is a subquasimodule $Q_2$ of $Q_1$
such that $Q_3=Q_1/Q_2$ is subdirectly irreducible and
non--associative. Now, $\gen(Q_3)=3$ and $Q_3$ is cocyclic by
7.4.
\qed\enddemo

\proclaim{7.7 Theorem} A finite quasimodule $Q$ is
non--associative and cocyclic if and only if there exists a
ternary algebra $A$ such that $Q=\q(A)$ {\rm(see 6.1)},
$\overline{\tau}\ne 0$ and the underlying module $A'=A(+,rx)$ is
cocyclic (in this case, $A'$ is $\widehat{\K}$--torsion and $\gen(A')\ge
3$).
\endproclaim

\demo{Proof} Assume first that $Q$ is both non--associative and
cocyclic. By 7.1 and 7.5, $Q$ is $\widehat{\K}$--torsion and
nilpotent of class 2. By 6.3, $Q$ has a ternary representation
$Q=\q(A)$. Now, the quasimodule $\q(A)$ and the module $A'$ have
the same cyclic submodules and it follows easily that $A'$ is
cocyclic and $\widehat{\K}$--torsion,. Further, by 6.1, we have
$\overline{\tau}\ne 0$ and one may check easily that then
$\gen(A')\ge 3$.

Now, the converse implication. Again, since $\q(A)$ and $A'$ have
the same cyclic submodules, $\q(A)$ is cocyclic. Finally, since
$\overline{\tau}\ne 0$, the quasimodule $\q(A)$ is
non--associative.
\qed\enddemo

\proclaim{7.8 Theorem} There exists a non--associative
hamiltonian quasimodule if and only if there exists a finite
cocyclic $\widehat{\K}$--torsion module $M$ such that $\gen(M)\ge
3$.
\endproclaim

\demo{Proof} The direct implication follows from 7.6 and 7.7 and
we have to show the converse one. To that purpose, we may assume
that $\gen(M)=3$. Further, consider the ternary algebra $F=F_3$
constructed in 6.2 and the corresponding free quasimodule $\q(F)$.
Let $D$ be a submodule of $_{\R}E=\R a_1+\R a_2+\R
a_3=\R^{(3)}\subseteq F$ such that $E/B\simeq\,_{\R}M$. Since $M$ is
$\widehat{\K}$--torsion, we have $\J(M)=(\I E+B)/B$, and
therefore $M/\J(M)\simeq E/(\I E+B)$. Since $M$ is finite and
$\gen(M)=3$, we have $\gen(M/\J(M))=3$, and consequently $E/(\I
E+B)$ is a $\K$-torsion module of dimension 3 and, in particular,
$|E/(\I E+B)|=27$. On the other hand, $\I E=\I a_1+\I a_2+\I a_3$
and $|E/\I E|=27$, too. Thus $B\subseteq \I E$ and, since $E/\I
E$ is not cocyclic, we have $B\ne\I E$ and $A/B\simeq\P$ for a
submodule $A$ of $\I E$, $B\subseteq A\subseteq \I E$. Now, fix
an epimorphism $\varphi:A\to\P$, $\Ker(\varphi)=B$, and define a
subset $V$ of $F$ by $x=(x_1,x_2,x_3,x_4)\in V$ if and only
if $f(x)=(x_1,x_2,x_3,0)\in A$ and $\varphi(f(x))=x_4$. Then
$V$ is a submodule of $\I E+\R a_4$ and, since $\I a_4=0$, we
have $u\*v=u+v$ for all $u,v\in V$. Consequently, $V$ is a
subquasimodule of $\q(F)$ and it is easy to check that $V$ is a
normal subquasimodule. We denote by $Q$ the corresponding
factorquasimodule; clearly $|Q|=3|M|$. Since $V\cap\R a_4=0$, the
quasimodule $Q$ is not associative. Furthermore, $P=(\R
a_4+V)/V\simeq\P$ is a normal simple submodule of $Q$ and
$Q/P\simeq \q(F)/(\R a_4+V)$ is a module. Consequently, $P=\A(Q)$
and, in order to show that $Q$ is hamiltonian, it is sufficient
to check that $P$ is contained in every non--zero cyclic
submodule of $Q$ or, equivalently, that $a_4\in\R x+V$ for every
$x\in F\setminus V$.

Let $x=(x_1,x_2,x_3,x_4)\in F\setminus V$ and
$y=f(x)=(x_1,x_2,x_3,0)$. Since $x\notin V$, we have either $y\in
A$ and $\varphi(y)\ne x_4$, or $y\notin A$. In the former case,
$z=(x_1,x_2,x_3,\varphi(y))\in V$, $0\ne x-z\in\R a_4$,
$a_4=r(x-z)$ for some $r\in\R$ and $a_4\in\R x+V$.

Assume that $y\notin A$. Since $M$ is finite and
$\widehat{\K}$--torsion, there is $m\ge 1$ with $\I^my\subseteq
B$ and $\I^{m-1}y\nsubseteq B$. If $m=1$ then $y+B\in\K(M)=A/B$
and $y\in A$, a~contradiction. Hence $m\ge 2$ and $sy\notin B$
for some $s\in\I^{m-1}$. Now, $sy\in A$ and
$sx=(sx_1,sx_2,sx_3,0)$. Finally,
$z=(sx_1,sx_2,sx_3,\varphi(sy))\in V$ and $z=sx+v$,
$v=(0,0,0,\varphi(sy))\in\R a_4$. Thus $v\in\R x+V$ and, since
$v\ne 0$, we conclude that $a_4\in\R x+V$.
\qed\enddemo

\flushpar{\bf 7.9} {\smc Remark.} If $\R$ is a principal ideal
domain then all finitely generated cocyclic modules are cyclic,
and hence every hamiltonian ($\R$--)quasimodule is a module.
\par\medskip

\flushpar{\bf 7.10} {\smc Remark.} According to \cite{20}, \cite{25}
and \cite{31}, there exists a finite (non--asso\-ciati\-ve) subdirectly
irreducible primitive quasimodule of order $n\ge 1$ if and only
if $n=3^m$ for $m\ge 1$, $m\ne 2,3,5$ ($m\ge 4$, $m\ne 5$).

\head 8. Trimedial quasigroups\endhead

Recall that by a trimedial quasigroup we mean a quasigroup $Q$
such that every subquasigroup $P$ generated by at most three
elements is medial (i.e., $P$ satisfies $(ax)(yb)=(ay)(xb)$
identically). We denote by $\Cal T$ the variety (or equational
class) of trimedial quasigroups and by $\Cal T^p$ that of pointed
trimedial quasigroups ($\Cal T^p$ contains just ordered pairs
$(Q,a)$, where $Q\in\Cal T$ and $a\in Q$).

The following basic result is proven in \cite{22} (see also
\cite{19}).

\proclaim{8.1 Proposition} The following conditions are
equivalent for a quasigroup $Q$:
\roster
\item"{\rm(i)}" $Q\in\Cal T$.
\item"{\rm(ii)}" There exist a commutative Moufang loop $Q(+)$
(defined on the same underlying set as $Q$), commuting 1-central
automorphisms $f,g$ of $Q(+)$ and a central element
$a\in\Z(Q(+))$ such that $xy=f(x)+g(y)+a$ for all $x,y\in Q$.\qed
\endroster
In this case, $Q$ is medial iff $Q(+)$ is associative.
\endproclaim

The ordered quadruple $(Q(+),f,g,a)$ will be called an {\sl
arithmetical form} of the trimedial quasigroup $Q$. Notice also
that $Q$ is medial if and only if $Q(+)$ is an abelian group.

\proclaim{8.2 Lemma} {\rm(\cite{22}, 3.2, 3.3)} Let $Q\in\Cal T$.
Then:
\roster
\widestnumber\item{}
\item"{\rm(i)}" For every $w\in Q$ there exists an arithmetical
form $(Q(+),f,g,a)$ of $Q$ such that $w=0$ is the neutral element
of the loop $Q(+)$ (then $a=ww$).
\item"{\rm(ii)}" If $(Q(+),f,g,a)$ and $Q(*),p,q,b)$ are
arithmetical forms of $Q$ such that the loops $Q(+)$ and $Q(*)$
possess the same neutral element then $Q(+)=Q(*)$, $f=p$, $g=q$
and $a=b$.\qed
\endroster
\endproclaim

\proclaim{8.3 Lemma} {\rm(\cite{22}, 3.4)} Let $Q(+),f,g,a)$ and
$P(+),p,q,b)$ be arithmetical forms of trimedial quasigroups $Q$
and $P$, respectively, and let $\varphi:Q\to P$ be a mapping such
that $\varphi(0)=0$. Then $\varphi$ is a homomorphism of the
quasigroups if and only if $\varphi$ is a homomorphism of the
loops such that $\varphi f=p\varphi$, $\varphi g=q\varphi$ and
$\varphi(a)=b$.\qed
\endproclaim

Put $\R_1=\Bbb Z[\x,\y,\x^{-1},\y^{-1}]$, $\x$ and $\y$ being two
commuting indeterminates over the ring $\Bbb Z$ of integers. Then
$\R_1$ is a commutative noetherian domain, a unique factorization
domain, and there exists just
one homomorphism $\F$ of $\R_1$ onto $\Bbb Z_3$; we have
$\F(\x)=2=\F(\y)$ and $\I=\Ker(\F)=3\R_1+(1+\x)\R_1+(1+\y)\R_1$.
Further, we denote by $\Cal Q_1^c$ the variety of centrally
pointed $\R_1$--quasimodules. That is, $\Cal Q_1^c$ contains just
ordered pairs $(\overline{Q},a)$, where $\overline{Q}$ is an
$\R_1$--quasimodule and $a\in\Z(\overline{Q})$.

\proclaim{8.4 Proposition} The varieties $\Cal T^p$ of pointed
trimedial quasigroups and $\Cal Q_1^c$ of centrally pointed
$\R_1$--quasimodules are equivalent.
\endproclaim

\demo{Proof} Let $(Q,w)\in\Cal T^p$. By 8.2(i), there is an
arithmetical form $(Q(+),f,g,a)$ of $Q$ such that $w=0$ is the
neutral element of $Q(+)$ and $a=ww\in\Z(Q(+))$. The
automorphisms $f,g$ are 1--central, i.e., $x+f(x)\in\Z(Q(+))$ and
$x+g(x)\in\Z(Q(+))$ for every $x\in Q$. Furthermore,
$3Q\subseteq\Z(Q(+))$ (which is true in every commutative Moufang
loop), and consequently we may turn $Q(+)$ into a quasimodule
$\overline{Q}$ by setting $\x x=f(x)$ and $\y x=g(x)$ for every
$x\in Q$; clearly, $\I Q\subseteq\Z(\overline{Q})$. Now,
$\lambda(Q,w)=(\overline{Q},a)\in\Cal Q_1^c$.

Conversely, take $(\overline{Q},a)\in\Cal Q_1^c$ and define a
binary operation on $Q$ by $xy=\x x+\y y+a$ for all $x,y\in Q$.
By 8.1, $Q$ becomes a trimedial quasigroup and we have
$\varkappa(\overline{Q},a)=(Q,0)\in\Cal T^p$.

We get correspondences $\lambda:\Cal T^p\to\Cal Q_1^c$,
$\varkappa:\Cal Q_1^c\to\Cal T^p$ and it follows easily from 8.2
that $\varkappa\lambda=\id$ and $\lambda\varkappa=\id$. Both
correspondences are biunique, they preserve the underlying sets
and, in view of 8.3, they represent equivalences between the
varieties.
\qed\enddemo

\proclaim{8.5 Lemma} Let $\alpha=(Q(+),f,g,a)$ and
$\beta=(Q(*),p,q,b)$ be arithmetical forms of a trimedial
quasigroup $Q$ (the neutral elements of $Q(+)$ and $Q(*)$ being
denoted by $0$ and $o$, respectively) and let $(\overline{Q},a)$
and $(\widetilde{Q},b)$ be centrally pointed quasimodules
corresponding to $\alpha$ and $\beta$, respectively {\rm(see
8.3)}. Then the quasimodules $\overline{Q}$ and $\widetilde{Q}$
are isomorphic.
\endproclaim

\demo{Proof} Define an operation $\circ$ on $Q$ by $x\circ
y=(x+y)-o$ for all $x,y\in Q$. Then $Q(\circ)$ is a loop, $o$ is
its neutral element and $h:Q(\circ)\to Q(+)$ is an isomorphism,
where $h(x)=x-o$. Moreover, $p_1=h^{-1}fh$ and $q_1=h^{-1}gh$ are
1-central automorphisms of $Q(\circ)$, $p_1q_1=q_1p_1$,
$b_1=f(o)+g(o)+a\in\Z(Q(\circ))$ and $xy=p_1(x)+q_1(y)+b_1$ for
all $x,y\in Q$. Now, by 8.2(ii), $Q(\circ)=Q(*)$, $p_1=p$,
$q_1=q$ and $b_1=b$, and hence $h$ is an isomorphism of the
quasimodules.
\qed\enddemo

\proclaim{8.6 Proposition} Let $Q$ be a trimedial quasigroup and
let $\overline{Q}$ be the corresponding quasimodule {\rm(see
8.4 and 8.5)}. Then:
\roster
\widestnumber\item{}
\item"{\rm(i)}" If $\overline{Q}$ is hamiltonian then $Q$ is so.
\item"{\rm(ii)}" If $Q$ is hamiltonian and contains at least one
idempotent element then $\overline{Q}$ is hamiltonian.
\item"{\rm(iii)}" $Q$ is medial iff $\overline{Q}$ is a module.
\endroster
\endproclaim

\demo{Proof} (i) Let $w\in P$, $P$ being a given subquasigroup of
$Q$, and let $(\widetilde{Q},a)$ be the centrally pointed
quasimodule corresponding to the pair $(Q,w)$ in the sense of 8.4.
Then $P$ is a subquasimodule of $\widetilde{Q}$ and, since
$\widetilde{Q}$ is hamiltonian, $P$ is a block of a~congruence
$r$ of $\widetilde{Q}$. Now, it is easy to check that $r$ is also
a normal congruence of the quasigroup $Q$.\newline
(ii) Let $e\in Q$ be such that $ee=e$ and let $(\widehat{Q},a)$
be the centrally pointed quasimodule corresponding to the pointed
quasigroup $(Q,e)$. Then $a=ee=e=0$, and so $e\in P$ for every
subquasimodule $P$ of $\widehat{Q}$. Now, $P$ is a normal
subquasigroup of $Q$ and the corresponding normal congruence $r$
of $Q$ is also a congruence of the quasimodule $\widehat{Q}$.
\qed\enddemo

\proclaim{8.7 Lemma} Let $Q$ be a trimedial quasigroup such that
the corresponding quasimodule $\overline{Q}$ {\rm(see 8.4, 8.5 and
5.5)} is subdirectly irreducible and nilpotent of class at most
2. Then every non--idempotent subquasigroup $P$ of $Q$ is a
normal subquasigroup.
\endproclaim

\demo{Proof} Take $w\in P$ such that $a=ww\ne w$ and let
$(\overline{Q},a)$, $\overline{Q}=Q(+,rx)$, be the centrally
pointed quasimodule corresponding to $(Q,w)$ in the sense of 8.4.
Clearly, $P(+)$ is a subloop of $Q(+)$ and $0=w\ne a\in
V=\Z(Q(+))\cap P$. Now, $V$ is a non--zero normal subquasimodule
of $\overline{Q}$, and hence $\A(Q(+))\subseteq V$. Thus
$\A(Q(+))\subseteq P$ and $P$ is a normal subquasimodule of
$\overline{Q}$. From this it easily follows that $P$ is a~normal
subquasigroup of $Q$.
\qed\enddemo

\flushpar{\bf 8.8} {\smc Remark.} The smallest possible number of
elements of a non--medial trimedial quasigroup is 81. According
to \cite{7}, there exist just 35 isomorphism classes of
non--medial trimedial quasigroups of order 81. Now, if $Q$ is
such a quasigroup and if $Q$ has no idempotent element
then $Q$ is hamiltonian by 8.7.
\par\medskip

\flushpar{\bf 8.9} {\smc Example.} Define an operation
$\diamondsuit$ on $\Cal D_1$ (see 3.1) by $a\diamondsuit
b=(a\bigtriangleup b)+(1,0,0,0)$. Then $\Cal D_1(\diamondsuit)$ is a
non-medial $CH$--quasigroup and $a\diamondsuit a=a+(1,0,0,0)\ne a$
for every $a\in\Cal D_1$. Thus $\Cal D_1(\diamondsuit)$ has no
idempotents and is hamiltonian (see 8.7, 8.8).
\par\medskip

\flushpar{\bf 8.10} {\smc Remark.} (i) In this remark, let us
call a quasigroup $Q$ {\sl meagre} ({\sl minimal}, resp.) if $Q$
is non--trivial and has no proper (non--trivial proper, resp.)
subquasi\-group.\newline
(ii) Every simple hamiltonian quasigroup is minimal. Conversely,
if $Q$ is minimal then $Q$ is hamiltonian and, moreover, if $Q$
contains at least one idempotent then $Q$ is simple.\newline
(iii) Every minimal trimedial quasigroup $Q$ is medial and,
moreover, $Q$ is either idempotent or contains just one
idempotent element or is meagre.\newline
(iv) Every simple trimedial quasigroup is minimal, finite and
medial (\cite{18} and \cite{19}).\newline
(v) Let $Q$ be a finite meagre quasigroup, $|Q|=q$, and $P$ be a
finite quasigroup such that $|P|=p$ is prime. Assume further that
the product $R=P\times Q$ is a hamiltonian quasigroup and
$\Hom(P,Q)=\emptyset$ (e.g., $P$ is meagre and not an image of
$Q$ or $P$ is meagre and $p$ does not divide $q$ or $P$ contains
no idempotent, $Q$ is simple and not isomorphic to a
subquasigroup of $P$). Then $R$ is not simple and we claim that
$R$ is meagre.

Indeed, if $S$ is a subquasigroup of $R$ then $s=|S|\ge q$ (since
$Q$ is meagre) and $s$ divides $|R|=qp$ (since $R$ is
hamiltonian). If $s>q$ then $s=qp$ (since $p$ is prime) and
$S=R$. On the other hand, if $s=q$ then, for every $a\in Q$, there
exists a unique $f(a)\in P$ with $(a,f(a))\in S$. Now, $f:Q\to P$
is a homomorphism, a contradiction.\newline
(vi) Put $R=\Bbb Z_5\times\Bbb Z_3$ and define an operation
$\circ$ on $R$ by $(a,x)\circ(b,y)=(3a+3b+1,2x+2y+1)$. Then
$R(\circ)$ is a commutative medial quasigroup and $R(\circ)$ is
meagre but not simple (see (iv)).

\head 9. Distributive quasigroups\endhead

Recall that a distributive quasigroup is characterized by the
equations $x(ab)=(xa)(xb)$ and $(ab)x=(ax)(bx)$ and that every
distributive quasigroup is trimedial (\cite{2}). Thus
distributive quasigroups are just idempotent trimedial
quasigroups.

Put $\R_2=\Bbb Z[\x,\x^{-1},(1-\x)^{-1}]$, $\x$ being an
indeterminate over the ring $\Bbb Z$ of integers. Then $\R_2$ is a
commutative noetherian domain and there exists just one
homomorphism $\F$ of $\R_2$ onto $\Bbb Z_3$; obviously, we have $\F(\x)=2$
and $\I=\Ker(\F)=3\R_2+(1+\x)\R_2$. Further, we denote by $\Cal
Q_2$ the variety of $\R_2$--quasimodules and by $\Cal D^p$ the
variety of pointed distributive quasigroups.

\proclaim{9.1 Proposition} The varieties $\Cal D^p$ of pointed
distributive quasigroups and $\Cal Q_2$ of~$\R_2$--quasimodules
are equivalent.
\endproclaim

\demo{Proof} Let $(Q,w)\in\Cal D^p$. By 8.2(i), there is an
arithmetical form $(Q(+),f,g,0)$ of $Q$ such that $w=0$ is the
neutral element of $Q(+)$, $f,g$ are 1--central automorphisms
of~$Q(+)$ and $x=xx=f(x)+g(x)$, $g(x)=(1-f)(x)$ for every $x\in Q$.
Consequently, we may turn $Q(+)$ into a quasimodule $\overline{Q}$
by setting $\x x=f(x)$ (see the proof of~8.4), and so
$\lambda(Q,w)=\overline{Q}\in\Cal Q_2$.

Conversely, if $\overline{Q}\in\Cal Q_2$ then
$(Q,0)=\varkappa(\overline{Q})$ is a pointed distributive
quasigroup, where the multiplication is defined by $xy=\x
x+(1\negmedspace-\negmedspace\x)y$ for all $x,y\in Q$.

Now, the correspondences $\lambda:\Cal D^p\to\Cal Q_2$ and
$\varkappa:\Cal Q_2\to\Cal D^p$ represent the desired equivalence
between the varieties (again, see the proof of 8.4).
\qed\enddemo

\proclaim{9.2 Proposition} Let $Q$ be a distributive quasigroup
and let $\overline{Q}$ be the corresponding quasimodule {\rm(see
9.1 and 8.5)}. Then:
\roster
\item"{\rm(i)}" $Q$ is hamiltonian if and only if
$\overline{Q}$ is so.
\item"{\rm(ii)}" $Q$ is medial iff $\overline{Q}$ is a module.
\endroster

\endproclaim

\demo{Proof} The assertion follows immediately from 8.6.
\qed\enddemo

\flushpar{\bf 9.3} {\smc Remark.} (i) Let $Q$ be a distributive
quasigroup and let $w_1,w_2\in Q$. Then $w_2=vw_1$ for some $v\in
Q$ and we have $w_2=\varphi(w_1)$, where $\varphi(x)=vx$ for
every $x\in Q$. Clearly, $\varphi$ is an automorphism of $Q$, and
so $\varphi$ is also an isomorphism of the pointed quasigroup
$(Q,w_1)$ onto the pointed quasigroup $(Q,w_2)$.\newline
(ii) There is a one--to--one correspondence between isomorphism
classes of distributive quasigroups and isomorphism classes of
$\R_2$--quasimodules. This correspondence preserves the
hamiltonian property.
\par\medskip

\flushpar{\bf 9.4} {\smc Remark.} Every non--medial distributive
quasigroup contains at least 81 elements and there exist just 6
isomorphism classes of non--medial distributive quasigroups of
order 81 (see \cite{27}).

\proclaim{9.5 Theorem} Let $A$ be a ternary $\R_2$--algebra such
that $\overline{\tau}\ne 0$ and the underlying module
$A'=A(+,rx)$ is cocyclic. Define an operation $\bigtriangledown$ on
$A$ by
$$
x\bigtriangledown y=\x
x+(1-\x)y+(\x^2+2\x^3)\tau(x,y,x)+(\x+\x^2+\x^3)\tau(y,x,y)
$$
for all $x,y\in A$. Then $A(\bigtriangledown)$ is a non--medial
hamiltonian distributive quasigroup.
\endproclaim

\demo{Proof} Combine 7.7, 9.1 and 9.2.
\qed\enddemo

\flushpar{\bf 9.6} {\smc Remark.} (i) A distributive quasigroup
$Q$ is simple if and only if $Q$ is non--trivial and contains no
non--trivial proper subquasigroup (i.e., $Q$ is minimal).\newline
(ii) Every simple distributive quasigroup is finite and medial.
\par\medskip

\flushpar{\bf 9.7} {\smc Remark.} (cf.\ \cite{11} and \cite{15})
Let $P$ be a minimal subquasigroup of a distributive quasigroup
$Q$.\newline
(i) Let $(Q(+),f,g,0)$ be an arithmetical form of $Q$ such that
$0\in P$. Then $P(+)$ is a minimal submodule of $Q(+)$ and $P$ is
a normal subquasigroup of $Q$ if and only if $P(+)$ is a normal
submodule of $Q(+)$ (use 9.1). Now, by 5.11, if $P$ is not normal
in $Q$ then $|P|=3$ and $|Pa\cdot Pb|\in\{3,9\}$ for all $a,b\in
Q$ (according to 2.3, we have $|Pa_0\cdot Pb_0|=9$ for some
$a_0,b_0\in Q$). Consequently (see 5.12), $P$ is normal in $Q$ if
and only if $|P|\cdot|Pa\cdot Pb|\ne 27$ for all $a,b\in
Q$.\newline
(ii) We show that $P\cdot ya=xy\cdot Pa$ for all $x,y\in P$ and
$a\in Q$. This is clear for $x=y$ and we assume that $x\ne y$.
Put $v=xy\cdot ya$ and $w=x\cdot ya=xy\cdot xa$. Then $v\ne w$,
$v,w\in P\cdot ya\cap xy\cdot Pa$ and $|P\cdot ya\cap xy\cdot
Pa|\ge 2$. But both $P\cdot ya$ and $xy\cdot Pa$ are minimal
subquasigroups of $Q$ and it follows that $P\cdot ya=xy\cdot
Pa$.\newline
(iii) It follows easily from (ii) that $P\cdot ya=P\cdot Pa$ for
all $y\in P$ and $a\in Q$. In particular, $|P\cdot Pa|=|P|$ and,
since $|p\cdot Pa|=|P|$ and $p\cdot Pa\subseteq P\cdot Pa$, we
have $p\cdot Pa=P\cdot Pa$ for every $p\in P$.\newline
(iv) Let $a,b\in Q$ be such that $Pa\cap Pb\ne\emptyset$. We show
that then $Pa=Pb$.
Indeed, $ua=q=vb$ for some $u,v\in P$ and, by (ii), $P\cdot
Pa=P\cdot ua=Pq=P\cdot vb=P\cdot Pb$, $p\cdot Pa=p\cdot Pb$ and,
finally, $Pa=Pb$.\newline
(v) According to (iv), $\{Pa\,|\,a\in Q\}$ is a partition of
$Q$.\newline
(vi) (\cite{15, 3.2}) Let $a,b\in Q$ be such that $|Pa\cdot
Pb|<|P|^2$. Then $xa\cdot yb=ua\cdot vb$ for some $x,y,u,v\in P$,
$(x,y)\ne(u,v)$, $P_1=Pa$ is a minimal subquasigroup of $Q$,
$P_1\cdot yb\cap P_1\cdot vb\ne\emptyset$, and hence $P_1\cdot
yb=P_1\cdot vb$ by (iv). If $c\in Q$ is such that $yb=xa\cdot c$
then $(pa)(P_1c)=P_1\cdot P_1c=P_1\cdot(xa\cdot c)=P_1\cdot
yb=P_1\cdot vb$ for every $p\in P$ (use (iii)). Now, it is clear
that $yb,vb\in Pb\cap P_1c$ and $|Pb\cap P_1c|\ge 2$.
Consequently, $Pb=P_1c=Pa\cdot c$ and $Pa\cdot Pb=Pa\cdot(Pa\cdot
c)=Pa\cdot(za\cdot c)$ for every $z\in P$ (again, by (iii)). Thus
$|Pa\cdot Pb|=|Pa\cdot(za\cdot c)|=|P|$.\newline
(vii) (cf.\ (i)) We have shown in (vi) that $|Pa\cdot Pb|$ is
equal to $|P|$ or $|P|^2$ for all $a,b\in Q$.
\par\medskip

\flushpar{\bf 9.8} {\smc Remark.} The variety of pointed
commutative distributive quasigroups is equivalent to the variety
of $\Bbb Z[\frac12]$-quasimodules. Since $\Bbb Z[\frac12]$ is a
principal ideal domain, every hamiltonian commutative
distributive quasigroup is medial (7.9).
\par\medskip

\flushpar{\bf 9.9} {\smc Remark.} (The parastrophes) Let $Q$ be an
$\R_2$-quasimodule. Keeping the underlying commutative Moufang
loop $Q(+)$ of $Q$, we introduce three new scalar
multiplications, say $\circ$, $\ast$ and $\bullet$, on $Q$ by the
equalities $\x\circ a=(1-\x)\cdot a$, $\x\ast a=\x^{-1}\cdot a$
and $\x\bullet a=-\x(1-\x)^{-1}\cdot a$. Then $(1-\x)\circ
a=\x\cdot a$, $(1-\x)\ast a=(1-\x^{-1})\cdot a$ and
$(1-\x)\bullet a=(1-x)^{-1}\cdot a$ and the resulting
quasimodules will be denoted by $\und\al(Q)=Q(+,r\circ a)$,
$\und\be(Q)=Q(+,r\ast a)$ and $\und\ga(Q)=q(+,r\bullet a)$,
respectively (the quasimodule $\und\al(Q)$ is called {\sl the
opposite quasimodule to $Q$} and is denoted also by
$\overline{Q}$). One may check easily that
$\und\al^2=\und\be^2=\und\ga^2=\id$,
$\und\al\und\be=\und\ga\und\al=\und\be\und\ga$ and
$\und\be\und\al=\und\al\und\ga=\und\ga\und\be$
(the six equivalences $\id$, $\und\al$, $\und\be$,
$\und\ga$, $\und\al\und\be$ and $\und\be\und\al$ form
a six-element group).

\head 10. Modules of divided powers\endhead

Throughout this section, let $\SS=\R[\x]$ and $\T=\R[[\x]]$
denote the ring of polynomials and the ring of formal power series
in one indeterminate $\x$ over $\R$, respectively. Now, given a
(unitary left) $\R$--module $M$, the direct sum
$N=M^{(\omega)}=\{a:\omega\to M\,|\,a(n)=0\text{ for almost all
}n\in\omega\}$ of $\omega$ copies of $M$ becomes a $\T$--module via
$(fa)(n)=\sum_{i=0}^{\infty}f_ia(n+i)$ for all $a\in N$ and
$f=\sum_{i=0}^{\infty}f_ix^i\in\T$ (the $\T$--module $_{\T}N$ is
known as the {\sl module of divided powers} and is denoted usually
by $M[\x^{-1}]$\,).

For $n\ge 0$, let $N_n=\{a\in N\,|\,a(m)=0\text{ for }m\ge
n+1\}$. Clearly, $N_0\subseteq N_1\subseteq N_2\subseteq\dots$ and
$N_n$ are submodules of all the three modules $_{\R}N$, $_{\SS}N$
and $_{\T}N$.

\proclaim{10.1 Lemma} $\SS a=\T a$ for every $a\in N$. If $a\ne 0$
then $\SS a\cap N_0\ne 0$.
\endproclaim

\demo{Proof} The equality is clear and if $n=\max\{\,i\,|\,a(i)\ne
0\}$ then $0\ne\x^na\in N_0$.
\qed\enddemo

\proclaim{10.2 Lemma} $\S(_{\T}N)=\S(_{\SS}N)=\S(_{\R}N_0)$.
\endproclaim

\demo{Proof} Every submodule of $_{\R}N_0$ is also a submodule of
$_{\T}N$, and hence $\S(_{\R}N_0)\subseteq\S(_{\T}N)$. On the
other hand, if $\T a$ is a (non--zero) simple submodule of
$_{\T}N$ then $\T a\cap N_0\ne 0$, and hence $\T a\subseteq N_0$.
Thus $\T a$ is a simple $\R$--module and $\T
a\subseteq\S(_{\R}N_0)$.
\qed\enddemo

\proclaim{10.3 Lemma} Let $m\ge 1$. Then
$\S_m(_{\T}N)=\S_m(_{\SS}N)\subseteq\S_m(_{\R}N_{m-1})\subseteq
\S_m(_{\R}N)$ and $a\in\S_m(_{\T}N)$ if and only if
$a(0)\in\S_m(_{\R}M)$, $a(1)\in\S_{m-1}(_{\R}M)$, \dots,
$a(m-1)\in\S_1(_{\R}M)$, $a(m)=a(m+1)=\dots=0$.
\endproclaim

\demo{Proof} We shall proceed by induction on $m$. The case $m=1$
is settled by 10.2, and so let $m\ge 2$.

First, take $b\in N$ such that $_{\T}B=(\T b+P)/P$ is a simple
$\T$--module, where $P=\S_{m-1}(_{\T}N)$. Then $_{\T}B\simeq\T/A$
for a maximal ideal $A$ of $\T$ and $\x^kb=0$ for some $k\ge 1$.
Consequently, since $\x^k\in A$ and $A$ is prime, we have $\x\in
A$, $\x B=0$ and $\x b\in P$. Now, by induction,
$b(1)\in\S_{m-1}(_{\R}M)$, \dots, $b(m-1)\in\S_1(_{\R}M)$,
$b(m)=b(m+1)=\dots=0$. Moreover, since $\x B=0$ and
$P\subseteq\S_{m-1}(_{\R}N)$, $B$ is also a simple $\R$--module,
$(\R b+\S_{m-1}(_{\R}N))/\S_{m-1}(_{\R}N)$ is a simple
$\R$--module, $b\in\S_m(_{\R}N)$ and $b(0)\in\S_m(_{\R}M)$.

Now, conversely, let $a\in N$ be such that
$a(i)\in\S_{m-i}(_{\R}M)$ for $0\le i\le m-1$ and $a(i)=0$ for
$i\ge m$. Then, by induction, $\x a\in P$, and hence $\T a+P=\R
a+P$. Moreover, if $C=(\R a+P)/P$ then $xC=0$ and $_{\R}C$ is
completely reducible. Consequently, $_{\T}C$ is also completely
reducible and $a\in\S_m(_{\T}N)$.
\qed\enddemo

\proclaim{10.4 Corollary}
$\S_{\omega}(_{\T}N)=\S_{\omega}(_{\SS}N)=\S_{\omega}(_{\R}N)
=\S_{\omega}(_{\R}M)^{(\omega)}$.\qed
\endproclaim

\proclaim{10.5 Lemma} If $m\ge 1$ and
$D_m=\S_m(_{\T}N)/\S_{m-1}(_{\T}N)$ then $\x D_m=0$ and
$_{\R}D_m\simeq\S_m(_{\R}M)/\S_{m-1}(_{\R}M)\times
\S_{m-1}(_{\R}M)/\S_{m-2}(_{\R}M)\times\dots\times
\S_1(_{\R}M)/\S_0(_{\R}M)$.
\endproclaim

\demo{Proof} The statements follow easily from 10.3.
\qed\enddemo

\proclaim{10.6 Lemma} Assume that
$\S_m(_{\R}M)\ne\S_{m-1}(_{\R}M)$ for some $m\ge 1$. Then none of
the modules $\S_m(_{\T}N)$ and $\S_m(_{\SS}N)$ can be generated by
less than $m$ elements.
\endproclaim

\demo{Proof} Since $\R$ is noetherian, $\widehat{\S}$--torsion
modules have primary decompositions, and hence there is a
homogeneous component $H$ of $\S_m(_{\R}M)$ such that
$H=\S_m(_{\R}H)\ne\S_{m-1}(_{\R}H)$. From this it follows that
$\S_{m-1}(_{\R}H)\ne\S_{m-2}(_{\R}H)$, \dots,
$\S_1(_{\R}H)\ne\S_0(_{\R}H)=0$ and consequently the module
$_{\R}D_m$ (see 10.5) contains a copy of $G^{(m)}$ for a simple
module $G$. The direct sum $G^{(m)}$ cannot be generated by $m-1$
elements and, since it is a direct summand of $_{\R}D_m$, the same
is true for the latter module. \qed\enddemo

\proclaim{10.7 Lemma} If $_{\R}M$ is cocyclic then both $_{\T}N$
and $_{\SS}N$ are cocyclic. Moreover, the modules $_{\R}M$,
$_{\T}N$ and $_{\SS}N$ are artinian and $\S_{\omega}$--torsion.
\endproclaim

\demo{Proof} See 10.1, \cite{44}, Theorem 4.30 and \cite{14}, Lemma
6.21.
\qed\enddemo

\proclaim{10.8 Lemma}{\rm(\cite{16}, Proposition~II.2)} The
following conditions are equivalent: \roster
\item"{\rm(i)}" $_{\R}M$ is an injective module.
\item"{\rm(ii)}" $_{\T}N$ is an injective module.
\item"{\rm(iii)}" ${}_{\SS}N$ is an injective module.
\endroster
\endproclaim

\demo{Proof} (i) $\Rightarrow$ (ii). This implication is
\cite{30}, Theorem 1.
\newline (ii) $\Rightarrow$ (iii). We will use a few
standard and well known arguments. First, since $\T$ is the
completion of $\SS$ in the usual $\x$--adic filtration, the module
${}_{\SS}T$ is flat (see e.g.\ \cite{36}, Theorem~8.8). Now, we
have the following natural transformation (see \cite{1},
Proposition~20.6):
$$
\Hom_{\SS}(_{\SS}A,\,_{\SS}N)\simeq
\Hom_{\SS}(_{\SS}A,\Hom_{\T}(_{\T}\T,\,_{\T}N))\simeq
\Hom_{\T}(_{\SS}T\otimes_{\SS}A,\,_{\T}N)
$$
for every $\SS$--module $A$. Since $_{\T}N$ is injective, the
functor $\Hom_{\T}(_{\SS}\T\otimes_{\SS}-,\,_{\T}N)$ is exact, and
hence the same is true for $\Hom_{\SS}(-,\,_{\SS}N)$. Thus
$_{\SS}N$ is injective.
\newline (iii) $\Rightarrow$ (i). We may
proceed in the same way as in the proof of \cite{30}, Prop.\ 1.
\qed\enddemo

\flushpar{\bf 10.9} Let $I$ be a maximal ideal of $\R$ and $E$ be
an injective envelope of the simple module $A=\R/I$. Then $K=\SS
x+\SS I$ is a maximal ideal of $\SS=\R[\x]$ and $B=\SS/K$ is a
simple $\SS$--module. Now, it follows from 10.7 and 10.8 that the
$\SS$--module $N=E[\x^{-1}]$ is an injective envelope of (a copy
of) $_{\SS}B$. The module $_{\SS}N$ is artinian,
$\S_{\omega}$--torsion and homogeneous, and every cocyclic
$\SS$--module containing $_{\SS}B$ (as the essential simple
socle) is isomorphic to a submodule of $_{\SS}N$. For every
$m\ge 1$, $\S_m(_{\SS}N)$ is a module of finite length (i.e.,
both artinian and noetherian) and if
$\S_m(_{\R}E)\ne\S_{m-1}(_{\R}E)$ (i.e., if
$E\ne\S_{m-1}(_{\R}E)$\,) then $\S_m(_{\SS}N)$ cannot be
generated by $m-1$ elements (use 10.6).
\par\medskip

\flushpar{\bf 10.10} Consider the situation from 10.9, take
$r_1\in\R$ and put $K_1=\SS(\x+r_1)+\SS I$ and $B_1=\SS/K_1$.
Clearly, $K_1$ is a maximal ideal of $\SS$ and $B_1$ is a simple
$\SS$--module. Now, denote by $\Cal B$ and $\Cal B_1$ the classes
of cocyclic $\SS$--modules $C$ and $C_1$, respectively, such that
$\S(C)\simeq B$ and $\S(C_1)\simeq B_1$. If $C\in\Cal B$ then
$\Lambda(C)=C_1\in\Cal B_1$, where both $\SS$--modules $C$ and
$C_1$ have the same underlying additive group and the
$\SS$--scalar multiplication $\cdot$ is defined on $C_1$ by $r\cdot
u=ru$ and $\x\cdot u=\x u-r_1u$ for all $r\in\R$ and $u\in C$.
Moreover, $\Lambda:\Cal B\to\Cal B_1$ is a bijective
correspondence and a subset $H$ of $C$ is a submodule of $C$ if
and only if it is a submodule of $C_1$. In particular, the
$\SS$--modules $C$ and $C_1$ possess the same number of
generators. Finally, a mapping $\varphi:C\to D$ is an
$\SS$--module homomorphism if and only if
$\varphi:\Lambda(C)\to\Lambda(D)$ is an $\SS$--module
homomorphism (it follows that $\Lambda$ is a category
equivalence).
\par\medskip

\flushpar{\bf 10.11} Let $C$ be a finite cocyclic $\SS$--module
such that $\x\S(C)\ne 0$. Then the mapping $u\mapsto \x u$ is an
automorphism of $C$, and hence $C$ becomes a (cocyclic)
$\R[\x,\x^{-1}]$--module. Similarly, if $\x v\ne v$ for at least
one $v\in\S(C)$, the mapping $u\mapsto(\x-1)u$ is an automorphism
of $C$ and $C$ is also an $\R[\x,(\x-1)^{-1}]$--module. Finally,
if $\x w\ne 0$ and $\x v\ne v$ for some $v,w\in\S(C)$ then $C$ is
an $\R[\x,\x^{-1},(1-\x)^{-1}]$--module.

\head 11. The socle series of $\Bbb
Z_{p^{\infty}}[\x^{-1}]$\endhead

This section is an immediate continuation of the preceding one.
Here, we choose $\R=\Bbb Z$, the ring of integers, and $\SS=\Bbb
Z[\x]$, the ring of polynomials with integral coefficients. For a
prime number $p\ge 2$, the module $N=\Bbb
Z_{p^{\infty}}[\x^{-1}]$ of divided powers is an injective
envelope of the simple $\SS$--module $B=\SS/(\SS\x+\SS p)$ (see
10.9) and $_{\SS}N$ is both artinian and $\S_{\omega}$--torsion.
Moreover, since $\S_m(_{\SS}N)\ne\S_{m-1}(_{\SS}N)$ for every
$m\ge 1$, the $m$--th member
$\S_m(_{\SS}N)$ of the socle series of $_{\SS}N$ cannot be
generated by $m-1$ elements; notice that
$|\S_m(_{\SS}N)|=p^{(m+1)m/2}$. Further, it is easy to see that
$\S_m(_{\SS}N)$ is isomorphic to the following $\SS$--module
$P_m$: $P_m=\Bbb Z_{p^m}\times\Bbb
Z_{p^{m-1}}\times\dots\times\Bbb Z_p$ and $(\x a)(n)=pa(n+1)$ for
$0\le n\le m-2$, $(\x a)(m-1)=0$, $a=(a(0),\dots,a(m-1))\in P_m$.
Clearly, the additive group $P_m(+)$ (and hence also the module
$_{\SS}P_m$) is generated by the elements $(1,0,\dots,0)$, \dots,
$(0,\dots,0,1)$ and $P_1\simeq\S/(\S3+\S\x)$..

\proclaim{11.1 Lemma}
$\J(\S_m(_{\SS}N))=\S_{m-1}(_{\SS}N)=p\S_m(_{\SS}N)$.
\endproclaim

\demo{Proof} The factor $\S_m(_{\SS}N)/\S_{m-1}(_{\SS}N)$ is a
completely reducible module isomorphic to $P_1^{(m)}$;
clearly, $\J(\S_m(_{\SS}N))\subseteq\S_{m-1}(_{\SS}N)$. On the
other hand, $\S_m(_{\SS}N)/\J(\S_m(_{\SS}N))$ is an
$m$--generated completely reducible module, and hence it is also
isomorphic to $P_1^{(m)}$. Thus
$\J(\S_m(_{\SS}N))=\S_{m-1}(_{\SS}N)$.
\qed\enddemo

\proclaim{11.2 Lemma} The $\SS$--module $\S_m(_{\SS}N)$ is
generated by $m$ elements but not by $m-1$ elements. Every proper
submodule of $\S_m(_{\SS}N)$ is generated by at most $m-1$
elements.
\endproclaim

\demo{Proof} Let $Q$ be a proper submodule of $P_m$. Assume first
that $Q+V=P_m$, where $V=\{a(0),0,\dots,0)\}\subseteq P_m$. Now, for
every $i$, $1\le i\le m-1$, there is $a_i(0)\in\Bbb Z_{p^m}$ such
that $a_1=(a_1(0),1,0,\dots,0)$, $a_2=(a_2(0),0,1,0,\dots,0)$,
\dots, $a_{m-1}=(a_{m-1}(0),0,\dots,0,1)$ are all in $Q$.
Clearly, $pV\subseteq Q_1\subseteq Q$, where $Q_1$ is the
submodule generated by the $m-1$ elements $a_1,\dots,a_{m-1}$. We
claim that $Q_1=Q$. Indeed, if $a\in Q$ then $a-b=c\in V$ for
some $b\in Q_1$. If $c\notin pV$ then $V=\SS c\subseteq Q$ and
$Q=P_m$, a contradiction. Thus $c\in pV\subseteq Q_1$ and $a\in
Q_1$.

Now, assume that $Q+V\ne P_m$. Then $P=Q/Q\cap V\simeq
Q+V/V\subseteq P_m/V\simeq P_{m-1}$, $P$ is isomorphic to a
proper submodule of $P_{m-1}$ and, using induction, we conclude
that $P$ is generated by at most $m-2$ elements. Finally, since
every submodule of $V$ is cyclic, $Q$ is generated by at most
$m-1$ elements.
\qed\enddemo

\proclaim{11.3 Proposition} Let $Q$ be a cocyclic $\SS$--module
whose (essential simple) socle is a copy of $P_1$ (i.e., $\Bbb
Z_p$, where $\x\Bbb Z_p=0$) such that
$\S_m(Q)=Q$. Then $Q$ is isomorphic to a submodule of
$\S_m(_{\SS}N)$ and $Q$ can be generated by at most $m$ elements.
Moreover, if $Q$ cannot be generated by $m-1$ elements then
$Q\simeq\S_m(_{\SS}N)$ {\rm($\simeq P_m$).}
\endproclaim

\demo{Proof} Since $Q$ is cocyclic and contains a copy of $P_1$,
$Q$ is isomorphic to a submodule of $_{\SS}N$. Further,
since $\S_m(Q)=Q$, a copy of $Q$ is contained in $\S_m(_{\SS}N)$
and the rest follows from 11.2.
\qed\enddemo

\proclaim{11.4 Lemma} Let $Q$ be a finitely generated cocyclic
$\SS$--module with $\S(Q)\simeq P_1$ and let $k=\gen(Q)$.
Then:
\roster
\widestnumber\item{}
\item"{\rm(i)}" If $k=1$ then $|Q|\ge p$, and if $|Q|=p$ then
$Q\simeq\S_1(_{\SS}N)\simeq P_1$.
\item"{\rm(ii)}" If $k=2$ then $|Q|\ge p^3$, and if $|Q|=p^3$
then $Q\simeq\S_2(_{\SS}N)\simeq P_2$.
\item"{\rm(iii)}" If $k=3$ then $|Q|\ge p^6$, and if $|Q|=p^6$
and $Q$ is not isomorphic to $\S_3(_{\SS}N)$ then the
$\S$--length of $Q$ is 4.
\item"{\rm(iv)}" If $k\ge 4$ then $|Q|\ge p^7$.
\endroster
\endproclaim

\demo{Proof} Easy (use 11.3; (iv) follows from (iii), since $Q$
contains a submodule $Q_1$ with $\gen(Q_1)=3$).
\qed\enddemo

\proclaim{11.5 Lemma}
$\S_1(P_4)$ is the set of all $a\in P_4$ such that $a(1)=a(2)=a(3)=0$
and $p^3$ divides $a(0)$, $\S_2(P_4)$ is the set of all $a\in P_4$
such that $a(2)=a(3)=0$ and $p^2$ divides $a(0),a(1)$, and
$\S_3(P_4)$ is the set of all $a\in P_4$ such that $a(3)=0$ and
$p$ divides $a(0),a(1),a(2)$.
\endproclaim

\demo{Proof} Easy.
\qed\enddemo

\proclaim{11.6 Lemma} Let $u,v\in\S_3(P_4)$ be any elements
such that the
submodule $Q=(\SS u+\SS v+\S(P_4))/\S(P_4)$ of $P_4/\S(P_4)$ is
not cyclic. Then
$\S(Q)=\S_1(P_4/\S(P_4))=\S_2(P_4)/\S_1(P_4)\simeq \Bbb Z_p^2$ is
not cyclic.
\endproclaim

\demo{Proof} By 11.5, we have $u=(i_0p,i_1p,i_2p,0)$ and
$v=(j_0p,j_1p,j_2p,0)$, where $0\le i_0,j_0\le p^3-1$, $0\le
i_1,j_1\le p^2-1$ and $0\le i_2,j_2\le p-1$. If at least one of
the elements $u,v,u-v$ is in $\S_2(P_4)$ then $|\S(Q)|=p^2$ (use
the fact that $Q$ is not cyclic), and hence we may assume that
none of these elements is in $\S_2(P_4)$ (see 11.5). Further, if
$i_1\ne 0\ne i_2$ then $(p^2,0,0,0)\in\SS u$, $(0,p^2,0,0)\in\SS
u$ and $\S_2(P_4)\subseteq\SS u$. Thus we may also assume that
$i_2=0$ implies $i_1=0$ and, similarly, $j_2=0$ implies $j_1=0$.
On the other hand, if $i_1=i_2=j_1=j_2=0$ then $\SS u+\SS v$ is a
cyclic module, a contradiction. Consequently, considering the
equalities $\SS u+\SS(u-v)=\SS u+\SS v=\SS v+\SS(u-v)$, we may
finally assume that $i_2=1$ and $j_1=j_2=0$. Then
$(p^2,0,0,0)\in\SS u+\SS v$, $(0,p^2,0,0)\in\SS u+\SS v $ and
$\S_2(P_4)\subseteq\SS u+\SS v$.
\qed\enddemo

\proclaim{11.7 Lemma} $P_2$ is not isomorphic to any submodule of
$P_4/\S_1(P_4)$.
\endproclaim

\demo{Proof} The result follows easily from 11.6.
\qed\enddemo

\proclaim{11.8 Lemma} Define an $\SS$--module structure on
$V=\Bbb Z_{p^2}\times\Bbb Z_{p^2}\times\Bbb Z_p$ by $\x
a=(pa(1),pa(2)),0)$. Then $P_2$ is not isomorphic to any
submodule of $_{\SS}V$.
\endproclaim

\demo{Proof} Let $u,v\in V$ be such that none of the elements
$u,v,u-v$ is in $\S(V)$ and $|\S(\SS u)|=p=|\S(\SS v)|$. Put
$Q=\SS u+\SS v$. One may check easily that either $|\S(Q)|=p^2$
or $Q$ is cyclic.
\qed\enddemo

\proclaim{11.9 Lemma} $P_2$ is not isomorphic to any submodule of
$P_4/\S_2(P_4)$.
\endproclaim

\demo{Proof} We have $P_4/\S_2(P_4)\simeq W=\Bbb
Z_{p^2}\times\Bbb Z_{p^2}\times\Bbb Z_{p^2}\times\Bbb Z_p$, where
the $\SS$--module structure is given by
$\x a=(pa(1),pa(2),pa(3),0)$. Now, let $Q$ be a submodule of
$W$ such that $Q\simeq P_2$. If $A=\{a\,|\,a(1)=a(2)=a(3)=0\}$
then $W/A\simeq V$ (see 11.8), and consequently $Q\cap A\ne 0$.
In particular, $\S(Q)=\SS(p,0,0,0)$. Finally, let $u,v\in Q$ be
such that $Q=\SS u+\SS v$. Then, using the fact that $Q$ is not
cyclic, we conclude easily that none of $u,v,u-v$ is in $\S(W)$
and $|\S(Q)|>p$, a contradiction.
\qed\enddemo

\proclaim{11.10 Proposition} Let $Q$ be a cocyclic $\SS$--module
with $\S(Q)\simeq P_1$ and $|Q|=p^6$. Then $Q$ is generated
by at most three elements, and if $\gen(Q)=3$ then
$Q\simeq P_3\simeq\S_3(_{\SS}N)$.
\endproclaim

\demo{Proof} Let $m$ denote the $\S$-length of $Q$. If $m\ge 2$
then (see 11.3) $Q$ is isomorphic to a submodule of
$\S_2(_{\SS}N)$ ($\simeq P_2$), and hence $|Q|\le p^3$, a
contradiction. Consequently, $m\ge 3$. Further, if $|\J(Q)|\le
p^2$ then $\J(Q)\subseteq\S_2(Q)$, $Q/\S_2(Q)$ is completely
reducible, $m=3$ and $\gen(Q)\le 3$ by 11.3 (a contradiction with
$|Q/\J(Q)|\ge p^4$). On the other hand, if $|\J(Q)|\ge p^3$ then
$|Q/\J(Q)|\le p^3$ and $\gen(Q)=\gen(Q/\J(Q))\le 3$. We have
proved that $\gen(Q)\le 3$. Now, assume that
$Q$ is not generated by two elements and that $Q$ is not
isomorphic to $P_3$. By 11.4(iii), $m=4$,
and hence, by 11.3, $Q$ is isomorphic to a submodule of $P_4$;
denote this submodule by $Q$ again and put $Q_1=Q/\S(P_4)$. Then
$|Q_1|=p^5$, the $\S$--length of $Q_1$ is 3 and, since
$\S(P_4)\subseteq\J(Q)$, the module $Q_1$ cannot be generated by
2 elements. Now, it follows from 11.3 that $Q_1$ is not cocyclic,
and therefore $|\S(Q_1)|\ge p^2$. On the other hand,
$\S(Q_1)\subseteq\S(P_4/\S(P_4))=\S_2(P_4)/\S_1(P_4)\simeq
P_1^{(2)}$ and we see that $|\S(Q_1)|=p^2$ and
$|Q_1/\S(Q_1)|=p^3$. Further, $Q_1/\J(Q_1)$ is a completely
reducible module which is 3--generated but not 2--generated,
$Q_1/\J(Q_1)\simeq P_1^{(3)}$, $|Q_1/\J(Q_1)|=p^3$ and
$|\J(Q_1)|=p^2$. If $\S(Q_1)\subseteq\J(Q_1)$ then
$\S(Q_1)=\J(Q_1)$, a contradiction with the fact that the
$\S$--length of $Q_1/\S(Q_1)$ is 2. Consequently,
$\S(Q_1)\nsubseteq\J(Q_1)$ and it follows that $Q_1=A\oplus Q_2$,
where $A,Q_2$ are submodules of $Q_1$ and $A$ is simple. Since
$\S(Q)\simeq P_1$, $Q_2$ is a cocyclic module and
$|Q_2|=p^4$. Clearly, $Q_2$ is 2--generated and not cyclic. The
$\S$--length of $Q_2$ is 3, and hence either $|\S_2(Q_2)|=p^3$ or
$|\S_2(Q_2)|=p^2$. If $|\S_2(Q_2)|=p^3$ then $\S_2(Q_2)$ is a
2--generated cocyclic module of $\S$--length 2 and
$\S_2(Q_2)\simeq P_2$ by 11.3. However this is a contradiction
with 11.7, and so $|\S_2(Q_2)|=p^2$ and $\S_2(Q_2)=\J(Q_2)$. Now,
$Q_2/\S(Q_2)\simeq P_2$ is isomorphic to a submodule of
$P_4/\S_2(P_4)$, which is, finally, a contradiction with 11.9.
\qed\enddemo

Define another scalar $\SS$--multiplication on $P_m$ by $(\x\cdot
a)(n)=pa(n+1)-a(n)$ for $0\le n\le m-2$, $(\x\cdot
a)(m-1)=-a(m-1)$, $a=(a(0),\dots,a(m-1))\in P_m$ (i.e., $\x\cdot
a=\x a-a=(\x-1)a$). In this way (see 10.10), we get a cocyclic
$\SS$--module $P_{m,1}$,
$S(P_{m,1})\simeq P_{1,1}\simeq\SS/(\SS p+\SS(\x+1))$.

\proclaim{11.11 Proposition} Let $Q_1$ be a finitely generated
cocyclic $\SS$--module such that $\S(Q_1)\simeq P_{1,1}$ and
$k=\gen(Q_1)$. Then:
\roster
\widestnumber\item{}
\item"{\rm(i)}" If $k\ge 3$ then $|Q_1|\ge p^6$.
\item"{\rm(ii)}" If $k\ge 4$ then $|Q_1|\ge p^7$.
\endroster
\endproclaim

\demo{Proof} Combine 11.4 and 10.10.
\qed\enddemo

\proclaim{11.12 Proposition} Let $Q_1$ be a cocyclic $\SS$--module
with $\S(Q_1)\simeq P_{1,1}$ and $|Q_1|=p^6$. Then $Q_1$ is
generated by at most three elements, and if $\gen(Q_1)=3$ then
$Q_1\simeq P_{3,1}$.
\endproclaim

\demo{Proof} Combine 11.10 and 10.10.
\qed\enddemo

\flushpar{\bf 11.13} {\smc Remark.} The transformation
$a\to\x\cdot a$ is an automorphism of $P_{m,1}(+)$ ($=P_m(+)$),
and if $p\ne 2$ then the same is true for the transformation
$a\to(\x-1)\cdot a$. Consequently, for $p\ne 2$, the scalar
$\SS$--multiplication on $P_{m,1}$ can be extended in a~unique way
to a scalar $\R_2$--multiplication (recall that $\R_2=\Bbb
Z[\x,\x^{-1},(1-\x)^{-1}]$) and the cocyclic $\SS$--module
$P_{m,1}$ turns into a cocyclic $\R_2$--module $P'_{m,1}$ (see
10.11). Notice that $P'_1\simeq \R_2/(\R_2p+\R_2(1+\x))$.

\proclaim{11.14 Proposition} Let $p\ne 2$ and let $Q'_1$ be a
finitely generated cocyclic $\R_2$--module with $S(Q'_1)\simeq
P'_{1,1}$ and $k=\gen(Q'_1)$. Then:
\roster
\widestnumber\item{}
\item"{\rm(i)}" If $k\ge 3$ then $|Q'_1|\ge p^6$.
\item"{\rm(ii)}" If $k\ge 4$ then $|Q'_1|\ge p^7$.
\endroster
\endproclaim

\demo{Proof} Combine 11.11 and 11.13.
\qed\enddemo

\proclaim{11.15 Proposition} Let $p\ne 2$ and let $Q'_1$ be a
cocyclic $\R_2$--module with $S(Q'_1)\simeq P'_{1,1}$ and
$|Q'_1|=p^6$. Then $Q'_1$ is generated by at most three elements,
and if $\gen(Q'_1)=3$ then $Q'_1\simeq P'_{3,1}$.
\endproclaim

\demo{Proof} Combine 11.12, 11.13 and use the fact that
$\x^{-1}\cdot a=\x^k\cdot a$ and $(1-\x)^{-1}\cdot
a=(1-\x)^l\cdot a$ for some positive integers $k$ and $l$.
\qed\enddemo

\head 12. The synthesis\endhead

\flushpar{\bf 12.1} In this (final) section, let $p=3$ and let
$\Cal P$ ($\Cal C$, resp.) denote the simple (cocyclic, resp.)
$\R_2$--module $P'_{1,1}$ ($P'_{3,1}$, resp.) defined in the
preceding section. Recall that $\Cal P(+)=\Bbb Z_3(+)$, $|\Cal
P|=3^1=3$, $3a=(1+\x)\cdot a=0$ and $\x\cdot a=\x^{-1}\cdot
a=(1-\x)\cdot a=(1-\x)^{-1}\cdot a=-a$ for every $a\in\Cal P$.
Further, $\Cal C(+)=\Bbb Z_{27}(+)\times\Bbb Z_9(+)\times\Bbb
Z_3(+)$, $|\Cal C|=3^6=729$, $3a=(3a(0),3a(1),0)$, $(1+\x)\cdot
a=(3a(1),3a(2),0)$, $\x\cdot a=(26a(0)+3a(1),8a(1)+3a(2),2a(2))$,
$\x^{-1}\cdot a=\x^{17}\cdot
a=(26a(0)+24a(1)+18a(2),8a(1)+6a(2),2a(2))$, $(1-\x)\cdot
a=(2a(0)+24a(1),2a(1)+6a(2),2a(2))$ and
$(1-\x)^{-1}\cdot a=(1-\x)^{17}\cdot
a=(14a(0)+21a(1)+18a(2),5a(1)+3a(2),2a(2))$ for every
$a=(a(0(,a(1),a(2))\in\Cal C$ (the transformations $a\to\x\cdot
a$ and $a\to(1-\x)\cdot a$ are permutations of $\Cal C$ and both
have order 18 in the corresponding symmetric group). We have
$\S(\Cal C)\simeq\Cal P$, $\gen(\Cal C)=3$ and, of course, $\Cal
C$ is $\widehat{K}$--torsion. By 11.15, $C\simeq\Cal C$ whenever
$C$ is a cocyclic $\R_2$--module with $\S(C)\simeq\Cal P$,
$\gen(C)=3$ and $|C|=729$. We put $\U=(1,0,0)$, $\V=(0,1,0)$ and
$\W=(0,0,1)$, $\U,\V,\W\in\Cal C$. Notice also that the mapping
$\lambda:a\mapsto(a(0),a(0)+8a(1),a(0)+a(1)+a(2))$ is an
automorphism of $\Cal C(+)$ such that $\lambda^2=\id$,
$\lambda(\x\cdot a)=(1-\x)\cdot\lambda(a)$ and
$\lambda((1-\x)\cdot a)=\x\cdot\lambda(a)$. Similarly, the
mapping $\ka:a\mapsto(a(0),-a(1)+6a(2),a(2))$ is an automorphism
of $\Cal C(+)$ such that $\ka(\x\cdot a)=\x^{-1}\cdot\ka(a)$,
$\ka((1-\x)\cdot a)=(1-\x^{-1})\cdot\ka(a)$ and the mapping
$\mu:a\mapsto(a(0)+9a(1),a(0)+a(1)+3a(2),a(0)+2a(1)+a(2))$ is an
automorphism of $\Cal C(+)$ such that $\mu(\x\cdot
a)=(1-\x)^{-1}\cdot\mu(a)$ and $\mu((1-\x)\cdot
a)=-\x(1-\x)^{-1}\cdot\mu(a)$. Finally, $\mu\la$ is an
automorphism of $\Cal C(+)$ such that $\mu\la(\x\cdot
a)=-\x(1-\x)^{-1}\cdot\mu\la(a)$ and $\mu\la((1-\x)\cdot
a)=(1-\x)^{-1}\cdot\mu\la(a)$.
\par\medskip

\flushpar{\bf 12.2} Let us define
$\tau_1(a,b,c)=(9a(0)b(1)c(2)+18a(1)b(0)c(2),0,0)$ and
$\tau_2(a,b,c)=(18a(0)b(1)c(2)+9a(1)b(0)c(2),0,0)$ for all
$a,b,c\in\Cal C$ (see 12.1). One checks readily that $\tau_1$,
$\tau_2$ are trilinear mappings of $\Cal C^{(3)}$ into $\Cal C$
and that these mappings satisfy the four conditions
T(0),\dots,T(3) from 1.5. Moreover, $\tau_2=-\tau_1=26\tau_1$,
$\overline{\tau_1}(a,b,c)=(9a(0)b(1)c(2)+9a(1)b(2)c(0)+
9a(2)b(0)c(1)+18a(0)b(2)c(1)+18a(1)b(0)c(2)+18a(2)b(1)c(0),0,0)
=-\overline{\tau_2}(a,b,c)$, $\overline{\tau_1}(\U,\V,\W)=(9,0,0)
=9\U$, $\overline{\tau_2}(\U,\V,\W)=(18,0,0)=18\U$. Thus $(\Cal
C,\tau_1)$ and $(\Cal C,\tau_2)$ are ternary algebras (see 1.5)
with $\overline{\tau_1}\ne 0\ne\overline{\tau_2}$.
\par\medskip

\flushpar{\bf 12.2.1} Define a mapping $\th:\Cal
C^{(3)}\mapsto\Cal C$ by
$\th(a,b,c)=(9(a(1)b(0)-a(0)b(1))(c(0)+c(1)),0,0)$ for all
$a,b,c\in\Cal C$. Clearly, $\th$ is a trilinear mapping
satisfying T(0),\dots,T(3) and $\overline{\th}=0$. Now,
consider the automorphism $\la$ of $\Cal C(+)$ defined in 12.1.
The following result is quite easy:

\proclaim{12.2.1.1 Lemma}
$\th(a,b,a-b)=(9(a(1)b(0)-a(0)b(1))(a(0)+a(1)-b(0)-b(1)),0,0)=
\tau_1(a,b,a-b)+\tau_1(\la(a),\la(b),\la(a-b))$ for all $a,b\in\Cal
C$.\qed\endproclaim

\flushpar{\bf 12.2.2} Put
$\si(a)=(18a(0)^3+9a(0)^2a(1)+18a(0)a(1)^2+9a(0),0,0)$ for
every $a\in\Cal C$.

\proclaim{12.2.2.1 Lemma}
$\si(a)+\si(b)=\si(a+b)+\th(a,b,a-b)$ for
all $a,b\in\Cal C$.
\endproclaim

\demo{Proof} Easy to check directly.
\qed\enddemo

\proclaim{12.2.2.2 Lemma} {\rm(i)} $\si(a)\ne 0$ if and
only if 3 divides $a(0)+a(1)$ and 3 does not divide $a(0)$ (or
$a(1)$).\newline
{\rm(ii)} If $\si(a)\ne 0$ then $\si(a)=9a$.
\endproclaim

\demo{Proof} Use the equality
$(a(0)+a(1)+1)(a(0)+a(1)+2)a(0)=a(0)^3+2a(0)^2a(1)+a(0)a(1)^2+
2a(0)+3a(0)^2+3a(0)a(1)$. \qed\enddemo

\proclaim{12.2.2.3 Corollary} $\si(a)\ne 0$ if and only if
either $a\in\{(1,2,0),(1,2,1),(1,2,2)\}+3\Cal C$ (and then
$\si(a)=(9,0,0)$) or
$a\in\{(2,1,0),(2,1,1),(2,1,2)\}+3\Cal C$ (and then
$\si(a)=(18,0,0)$).
\qed\endproclaim

\proclaim{12.2.2.4 Corollary} $|\{a\,|\,\si(a)=0\}|=567$ and
$|\{a\,|\,\si(a)\ne 0\}|=162$.
\qed\endproclaim

\flushpar{\bf 12.2.3} Put
$\xi(a)=\la(a)+\si(a)=(18a(0)^3+9a(0)^2a(1)+18a(0)a(1)^2
+10a(0),a(0)+8a(1),a(0)+a(1)+a(2))$ for every $a\in\Cal C$.

\proclaim{12.2.3.1 Lemma} $\xi^2(a)=a+2\si(a)$ for every
$a\in\Cal C$.
\endproclaim

\demo{Proof} We have $\xi^2(a)=\xi(\la(a)+\si(a))=
\la(\la(a)+\si(a))+\si(\la(a)+\si(a))=
\la^2(a)+\la\si(a)+\si\la(a)+\si^2(a)
-\th(\la(a),\si(a),\la(a)-\si(a))$ by 12.2.2.1. On the other
hand, $\la^2=\id$, $\la\si(a)=\si(a)=\si\la(a)$,
$\si^2(a)=0=\th(\la(a),\si(a),\la(a)-\si(a))$, and so
$\xi^2(a)=a+2\si(a)$.
\qed\enddemo

\proclaim{12.2.3.2 Corollary} $\xi^2(a)=a$ if and only if
$\si(a)=0$.
\qed\endproclaim

\proclaim{12.2.3.3 Lemma} $\xi$ is a permutation of $\Cal C$.
\endproclaim

\demo{Proof} Since $\Cal C$ is finite, it suffices to show that
$\xi$ is injective. However, if $\xi(a)=\xi(b)$ then
$a+2\si(a)=\xi^2(a)=\xi^2(b)=b+2\si(b)$, and hence $a(1)=b(1)$ and
$a(2)=b(2)$. Further, $3a=3a+6\si(a)=3b+6\si(b)=3b$, and so 27
divides $3(a(0)-b(0))$. From this, 27 divides
$18(a(0)^3-b(0)^3)$, etc., and we conclude that
$0=a+2\si(a)-b-2\si(b)=(19(a(0)-b(0)),0,0)$. Thus $a(0)=b(0)$ and
$a=b$.
\qed\enddemo

\proclaim{12.2.3.4 Lemma} $\xi(\x\cdot a)=(1-\x)\cdot\xi(a)$ and
$\xi((1-\x)\cdot a)=\x\cdot\xi(a)$ for every $a\in\Cal C$.
\endproclaim

\demo{Proof} We have $\xi(\x\cdot a)=\la(\x\cdot a)+\si(\x\cdot
a)=(1-\x)\cdot\la(a)+\si(\x\cdot a)$. It is easy to check that
$\si(\x\cdot a)=-\si(a)=(1-\x)\cdot\si(a)$, and hence
$\xi(\x\cdot
a)=(1-\x)\cdot\la(a)+(1-\x)\cdot\si(a)=(1-\x)\cdot\xi(a)$. Quite
similarly, $\xi((1-\x)\cdot a)=\x\cdot\xi(a)$.
\qed\enddemo

\proclaim{12.2.3.5 Lemma} $\xi(a+b+\tau_1(a,b,a-b))=\xi(a)+\xi(b)
+\tau_2(\xi(a),\xi(b),\xi(a)-\xi(b))$ for all $a,b\in\Cal C$.
\endproclaim

\demo{Proof} Using 12.2.1.1 and 12.2.2.1, one checks easily that
$\xi(a+b+\tau_1(a,b,a-b))=\la(a)+\la(b)+\tau_1(a,b,a-b)+\si(a+b)
=\xi(a)+\xi(b)+\tau_1(a,b,a-b)-\th(a,b,a-b)=\xi(a)+\xi(b)-\tau_1
(\la(a),\la(b),\la(a)-\la(b))=\xi(a)+\xi(b)-\tau_1(\xi(a),\xi(b),
\xi(a)-\xi(b))$.
\qed\enddemo

\flushpar{\bf 12.3} Define operations $\*$ and $\boxtimes$ on
$\Cal C$ (see 12.1, 12.2) by
$a\*b=a+b+\tau_1(a,b,a-b)=(a(0)+b(0)+9a(0)a(2)b(1)+18a(0)b(1)b(2)
+18a(1)a(2)b(0)+9a(1)b(0)b(2),a(1)+b(1),a(2)+b(2))$ and
$a\boxtimes
b=a+b+\tau_2(a,b,a-b)=(a(0)+b(0)+18a(0)a(2)b(1)+9a(0)b(1)b(2)
+9a(1)a(2)b(0)+18a(1)b(0)b(2),a(1)+b(1),a(2)+b(2))$ for all
$a,b\in\Cal C$. By 6.1, $\Cal C(\*)$ and $\Cal C(\boxtimes)$ are
commutative Moufang loops and it is easy to see that the mapping
$(a(0),a(1),a(2))\to(a(0),a(1),2a(2))$ is an isomorphism of $\Cal
C(\*)$ onto $\Cal C(\boxtimes)$. In fact, by \cite{27, 6.3},
these (isomorphic) loops
are determined by three generators, say $\alpha,\beta,\gamma$,
and relations $27\alpha=9\beta=3\gamma=0$,
$[\alpha,\beta,\gamma]=9\alpha$ (or
$[\alpha,\beta,\gamma]=18\alpha$). By 12.2.3.2 and 12.2.3.5,
$\xi:\Cal C(\*)\mapsto\Cal C(\boxtimes)$ is another iso\-mor\-phism
of the loops. Notice that $\xi(\U)=(1,1,1)$, $\xi(\V)=(0,-1,1)$,
$\xi(\W)=(0,0,1)$, $[\,\xi(\U),\xi(\V),\xi(\W)]_{\Cal C(\*)}=
\overline{\tau_1}((1,1,1),(0,-1,1),(0,0,1))=(18,0,0)=
18\U=18\xi(\U)$ and $[\,\xi(\U),\xi(\V),\xi(\W)]_{\Cal
C(\boxtimes)}=(9,0,0)=9\U=9\xi(\U)$.
\par\medskip

\flushpar{\bf 12.4} Put $\Cal C'_1=\q(\Cal C,\tau_1)=\Cal
C(\*,rx)$ and $\Cal C'_2=\q(\Cal C,\tau_2)=\Cal C(\boxtimes,rx)$
(see 6.1 and 12.3). By 7.7, both $\Cal C'_1$ and $\Cal C'_2$ are
non--associative cocyclic ($\widehat{K}$--torsion)
$\R_2$--quasimodules and we have $\gen(\Cal C'_1)=3=\gen(\Cal
C'_2)$.

\proclaim{12.4.1 Lemma}
\roster
\runinitem"{\rm(i)}" $27\U=9\V=3\W=0$.
\widestnumber\item{}
\item"{\rm(ii)}" $(1+\x)\cdot\U=0$, $(1+\x)\cdot\V=3\U$ and
$(1+\x)\cdot\W=3\V$.
\item"{\rm(iii)}" $[\U,\V,\W]_{\Cal C'_1}=9\U$ and
$[\U,\V,\W]_{\Cal C'_2}=18\U$.
\endroster
\endproclaim

\demo{Proof} See 12.1, 12.2 and 12.3.
\qed\enddemo

\proclaim{12.4.2 Lemma} The quasimodules $\Cal C'_1$ and $\Cal
C'_2$ are not isomorphic.
\endproclaim

\demo{Proof} Let, on the contrary, $\varphi:\Cal C'_2\to\Cal C'_1$
be an isomorphism. Put $u=\varphi(\U)$, $v=\varphi(\V)$ and
$w=\varphi(\W)$. We have $(3u(1),3u(2),0)=(1+\x)\cdot
u=\varphi((1+\x)\cdot\U)=0$, and hence $u(1)=0=u(2)$. Similarly,
$(3v(1),3v(2),0)=(1+\x)\cdot
v=\varphi((1+\x)\cdot\V)=\varphi(3\U)=3u=(3u(0),0,0)$, 9 divides
$u(0)-v(1)$, $v(2)=0$, $(3w(1),3w(2),0)=(1+\x)\cdot
w=\varphi((1+\x)\cdot\W)=\varphi(3\V)=3v=(3v(0),3v(1),0)$, 9
divides $v(0)-w(1)$, 3  divides $v(1)-w(2)$. Furthermore, the
orders of $u$, $v$ and $w$ are 27, 9 and 3, resp., and so 3 does
not divide $u(0)$, 3 divides $v(0)$, 3 does not divide $v(1)$, 9
divides $w(0)$, 3 divides $w(1)$ and $w(2)\ne 0$. Finally,
$(9u(0)v(1)w(2),0,0)=\overline{\tau_1}(u,v,w)=[u,v,w]_{\Cal
C'_1}=\varphi([\U,\V,\W]_{\Cal
C'_2})=\varphi(18\U)=18u=(18u(0),0,0)$, hence
$9u(0)(2-v(1)w(2))=0\text{ (mod 27)}$ and 3 divides $2-v(1)w(2)$.
Thus we have shown that 3 divides $2-w(2)^2$, since
$2-w(2)^2=2-w(1)w(2)+w(1)w(2)-w(2)^2$. But this is a contradiction
with $w(2)=1$ or 2. \qed\enddemo

\proclaim{12.4.3 Lemma} The permutation $\xi$ of $\Cal C$ {\rm(see
12.2.3)} is an anti--isomorphism of the quasimodule $\Cal C'_1$
onto the quasimodule $\Cal C'_2$.
\endproclaim

\demo{Proof} By 12.2.3.4 and 12.2.3.5, we have
$\xi(a\*b)=\xi(a)\boxtimes\xi(b)$, $\xi(\x\cdot
a)=(1-\x)\cdot\xi(a)$, $\xi((1-\x)\cdot a)=\x\cdot\xi(a)$ for all
$a,b\in\Cal C$ and this means that $\xi$ is an anti--isomorphism
of the quasimodules.
\qed\enddemo

\flushpar{\bf 12.4.4} {\smc Remark.} Let $\overline{\Cal C'_1}$
($\overline{\Cal C'_2)}$, resp.) denote the quasimodule opposite
to $\Cal C'_1$ ($\Cal C'_2$, resp.). That is, $\Cal C(\*)$ ($\Cal
C(\boxtimes)$, resp.) is the underlying commutative Moufang loop
of $\overline{\Cal C'_1}$ ($\overline{\Cal C'_2}$, resp.) and the
(opposite) scalar multiplication, say $\circ$, is given (in both
cases) by $\x\circ a=(1-\x)\cdot a$ and $(1-\x)\circ a=\x\cdot
a$. If $\U_1=(1,1,1)$ and $\V_1=(0,-1,1)$ then $\overline{\Cal
C'_1}$ ($\overline{\Cal C'_2}$, resp.) is generated by
$\{\U_1,\V_1,\W\}$, $27\U_1=9\V_1=3\W=0$, $(1+\x)\circ\U_1=0$,
$(1+\x)\circ\V_1=3\U_1$ and $(1+\x)\circ\W=3\V_1$, and
$[\U_1,\V_1,\W]_{\Cal C(\*)}=18\U_1$ ($[\U_1,\V_1,\W]_{\Cal
C(\boxtimes)}=9\U_1$, resp.). Now, it follows from 12.5.1 that
there exist quasimodule isomorphisms $\zeta_1:\Cal
C'_2\mapsto\overline{\Cal C'_1}$, $\zeta_2:\Cal C'_1\mapsto
\overline{\Cal C'_2}$ such that $\zeta_1(\U)=\U_1=\zeta_2(\U)$,
$\zeta_1(\V)=\V_1=\zeta_2(\V)$ and
$\zeta_1(\W)=\W=\zeta_2(\W)$. Then, of course, $\zeta_1$ is an
anti--isomorphism of $\Cal C'_2$ onto $\Cal C'_1$, $\zeta_2$ is
an anti--isomorphism of $\Cal C'_1$ onto $\Cal C'_2$ and, in fact,
$\zeta_1=\xi^{-1}=\zeta_2^{-1}$, $\zeta_2=\xi=\zeta_1^{-1}$
(use 12.2.3.2 and 12.4.3).
\par\medskip

\flushpar{\bf 12.4.5} {\smc Remark.} (i) Putting $\W_1=(0,6,1)$,
we have $27\U=9(-\V)=3\W_1=0$, $(1+\x^{-1})\cdot\U=0$,
$(1+\x^{-1})\cdot(-\V)=3\U$, $(1+\x^{-1})\cdot\W_1=-3\V$,
$[\U,-\V,\W_1]_{\Cal C(\*)}=18\U$ and $[\U,-\V,\W_1]_{\Cal
C(\boxtimes)}=9\U$. Consequently, considering the parastrophes
$\und\be(\Cal C'_1)$, $\und\be(\Cal C'_2)$ (see 9.8) and 12.5.1,
we get quasimodule isomorphisms $\eta_1:\Cal
C'_2\mapsto\und\be(\Cal C'_1)$ and $\eta_2:\Cal
C'_1\mapsto\und\be(\Cal C'_2)$ such that
$\eta_1(\U)=\U=\eta_2(\U)$, $\eta_1(\V)=-\V=\eta_2(\V)$ and
$\eta_1(\W)=\W_1=\eta_2(\W)$. In fact, $\eta_2=\eta_1^{-1}$.
\newline
(ii) Put $\U_2=(10,5,1)$, $\V_2=(18,2,2)$ and $\W_2=(0,3,1)$.
Then we have $27\U_2=9\V_2=3\W_2=0$,
$(1-2\x)(1-\x)^{-1}\cdot\U_2=0$ ($=(1-2\x)\cdot\U_2$),
$(1-\x)^{-1}\cdot\U_2=2\U_2$, $(1-2\x)(1-\x)^{-1}\cdot\V_2=3\U_2$,
$(1-2\x)(1-\x)^{-1}\cdot\W_2=3\V_2$, $[\U_2,\V_2,\W_2]_{\Cal
C(\*)}=18\U_2$ and $[\U_2,\V_2,\W_2]_{\Cal C(\boxtimes)}=9\U_2$.
Similarly as in (i), we get quasimodule isomorphisms
$\ro_1:\Cal C'_2\mapsto\und\ga(\Cal C'_1)$ and
$\ro_2:\Cal C'_1\mapsto\und\ga(\Cal C'_2)$ such that
$\ro_1(\U_)=\U_2=\ro_2(\U)$, $\ro_1(\V)=\V_2=\ro_2(\V)$ and
$\ro_1(\W)=\W_2=\ro_2(\W)$. Again, $\ro_2=\ro_1^{-1}$.
\newline
(iii) According to (i), (ii), 12.4.3, 12.4.4 and 9.9, we have the
following quasimodule isomorphisms: $\Cal C'_1\simeq\und\al(\Cal
C'_2)\simeq\und\be(\Cal C'_2)\simeq\und\ga(\Cal
C'_2)\simeq\und\al\und\be(\Cal C'_1)\simeq\und\be\und\al(\Cal
C'_1)\simeq\und\al\und\ga(\Cal C'_1)\simeq\und\ga\und\al(\Cal
C'_1)\simeq\und\be\und\ga(\Cal C'_1)\simeq\und\ga\und\be(\Cal
C'_1)$ and $\Cal C'_2\simeq\und\al(\Cal
C'_1)\simeq\und\be(\Cal C'_1)\simeq\und\ga(\Cal
C'_1)\simeq\und\al\und\be(\Cal C'_2)\simeq\und\be\und\al(\Cal
C'_2)\simeq\und\al\und\ga(\Cal C'_2)\simeq\und\ga\und\al(\Cal
C'_2)\simeq\und\be\und\ga(\Cal C'_2)\simeq\und\ga\und\be(\Cal
C'_2)$.
\par\medskip

\flushpar{\bf 12.5} Let $F$ be a free $\R_2$--quasimodule freely
generated by a three--element set $\{\alpha,\beta,\gamma\}$ (see
6.2). Then $F$ is nilpotent of class at most 2 and we denote by
$G_1$ ($G_2$, resp.) the subquasimodule generated by the elements
$27\alpha$, $9\beta$, $3\gamma$, $(1+\x)\alpha$,
$(1+\x)\beta\ominus3\alpha$, $(1+\x)\gamma\ominus3\beta$ and
$[\alpha,\beta,\gamma]\ominus9\alpha$
($[\alpha,\beta,\gamma]\ominus18\alpha$, resp.). Then
$G_1\subseteq\Z(F)$ ($G_2\subseteq\Z(F)$, resp.), and hence
$G_1$ ($G_2$, resp.) is a normal submodule of $F$.

\proclaim{12.5.1 Lemma} $F/G_1)\simeq\Cal C'_1$
($F/G_2\simeq\Cal C'_2$, resp.).
\endproclaim

\demo{Proof} The assertion follows easily from 6.2, but we may
also proceed in the following way: First, 12.4.1 implies that
$\Cal C'_1\simeq F/H$ for a normal subquasimodule $H$ of $F$ with
$G_1\subseteq H$. If $K$ is the (normal) subquasimodule of $F$
such that $G_1\subseteq K$ and $K/G_1)=\A(F/G_1)$ then
$[\alpha,\beta,\gamma]\in K$, $M=G_1/K$ is a module and if
$\pi:F\to M$ is the natural projection then
$9\pi(\alpha)=9\pi(\beta)=3\pi(\gamma)=0$,
$\x\pi(\alpha)=-\alpha$, $\x\pi(\beta)=3\alpha-\beta$,
$\x\pi(\gamma)=3\beta-\gamma$. It follows easily that the
additive group $M(+)$ is generated by the elements
$\pi(\alpha),\pi(\beta),\pi(\gamma)$ and $|M|\le3^5=243$. On the
other hand, $|K/G_1)|=|\A(F/G_1)|=3$, and therefore
$729\ge|F/G_1|\ge|F/H|=729$, $G_1=H$ and $F/G_1\simeq\Cal C'_1$.
\qed\enddemo

\proclaim{12.6 Theorem} Every non--associative cocyclic
$\R_2$--quasimodule contains at least 729 elements and the only
non--associative cocyclic $\R_2$--quasimodules of order 729 are (up
to isomorphism) the non--isomorphic quasimodules $\Cal C'_1$ and
$\Cal C'_2$ {\rm (see 12.4)}. These two quasimodules are
anti--isomorphic.
\endproclaim

\demo{Proof} The first assertion follows by easy combination of
7.7 and 11.11. Now, let $Q$ be a non--associative cocyclic
$\R_2$--quasimodule of order 729. By 7.7, there exists a ternary
algebra $A$ ($=A(+,rx,\tau)$) such that $Q=\q(A)$,
$\overline{\tau}\ne 0$, the module $A'=A(+,rx)$ is cocyclic and
$\widehat{K}$--torsion, and $\gen(A')\ge 3$. By 11.15, we have
$A'\simeq P'_{3,1}$ and we can assume that $A'=P'_{3,1}$ (see
11.13). We have $27\U=9\V=3\W=0$, $(1+\x)\cdot\U=0$,
$(1+\x)\cdot\V=3\U$, $(1+\x)\cdot\W=3\V$ and, since $Q$ is
cocyclic and $9Q$ is a non--trivial normal subquasimodule of $Q$,
we also have $\A(Q)\subseteq 9Q$ and $[\U,\V,\W]\in\{9\U,18\U\}$.
The subquasimodule $P$ of $Q$ generated by $\{\U,\V,\W\}$ is
non--associative and cocyclic, hence $|P|\ge 729$ and necessarily
$P=Q$. Now, it follows from 12.5.1 that $Q$ is a homomorphic
image of $\Cal C'_1$ or $\Cal C'_2$, and consequently either
$Q\simeq\Cal C'_1$ or $Q\simeq\Cal C'_2$. Finally, $\Cal C'_1$
and $\Cal C'_2$ are not isomorphic (12.4.2) but they are
anti--isomorphic (12.4.3; see also 12.4.4).
\qed\enddemo

\flushpar{\bf 12.7} Define two binary operations $\bigtriangledown_1$
and  $\bigtriangledown_2$ on $\Cal C$ by
$a\bigtriangledown_1 b=\x\cdot a\*(1-\x)\cdot b=
(26a(0)+3a(1),8a(1)+3a(2),2a(2))\*
(2b(0)+24b(1),2b(1)+6b(2),2b(2))=(26a(0)+3a(1)+2b(0)+24b(1)+
18a(0)a(2)b(1)+9a(0)b(1)b(2)+9a(1)a(2)b(0)+$ $18a(1)b(0)b(2),
8a(1)+3a(2)+2b(1)+6b(2),2a(2)+2b(2))$
and
$a\bigtriangledown_2 b=\x\cdot a\boxtimes(1-\x)\cdot b=
(26a(0)+3a(1),8a(1)+3a(2),2a(2))\boxtimes
(2b(0)+24b(1),2b(1)+6b(2),2b(2))=(26a(0)+3a(1)+2b(0)+24b(1)+
9a(0)a(2)b(1)+18a(0)b(1)b(2)+18a(1)a(2)b(0)+9a(1)b(0)b(2)$,$
8a(1)+3a(2)+2b(1)+6b(2),2a(2)+2b(2))$
for all $a,b\in\Cal C$; we have
$a\bigtriangledown_2 b=(a\bigtriangledown_1 b)-\tau_1(a,b,a-b)=
(a\bigtriangledown_1)+\tau_2(a,b,a-b)$ and $\Cal
C(\bigtriangledown_1)=\Cal D_2(\bigtriangledown)$ (see 3.2).
The permutation $\xi$ of
$\Cal C$ is an anti--isomorphism of $\Cal C(\bigtriangledown_1)$
onto $\Cal C(\bigtriangledown_2)$, i.e. $\xi(a\bigtriangledown_1
b)=\xi(b)\bigtriangledown_2\xi(a)$ for all $a,b\in\Cal C$ (see
12.2.3 and 12.4.3).

\proclaim{12.8 Theorem}
{\rm(i)}
Every non--medial hamiltonian distributive
quasigroup has at least 729 elements.
\newline{\rm(ii)}
$\Cal C(\bigtriangledown_1)$ and $\Cal
C(\bigtriangledown_2)$ are (up to isomorphism) the only non--medial
hamiltonian distributive quasigroups of order 729; these two
quasigroups are not isomorphic, but they are anti--isomorphic.
\endproclaim

\demo{Proof} Combine 7.7, 9.1, 9.2, 9.3, 12.6 and 12.7.
\qed\enddemo

\flushpar{\bf 12.9} {\smc Remark.} The preceding theorem says
that up to usual equivalences (as isomorphism and parastrophy)
there exists only one non-medial hamiltonian distributive
quasigroup of order 729, which is the smallest possible order for
such a~structure. The mapping
$a\mapsto(18a(0)^3+9a(0)^2a(1)+18a(0)a(1)^2+10a(0),a(0)+8a(1),
a(0)+a(1)+a(2))$ is an anti--isomorphism of $\Cal
C(\bigtriangledown_1)$ onto $\Cal C(\bigtriangledown_2)$, and so
$\Cal C(\bigtriangledown_2)$ is isomorphic to the opposite
quasigroup $\overline{\Cal C(\bigtriangledown_1)}$.
\par\medskip

\flushpar{\bf 12.10} {\smc Remark.} Using 12.4.5 (and also 12.4.4
and 12.9), we come to the following isomorphisms for the
parastrophes of the quasigroups $\Cal C(\bigtriangledown_1)$ and
$\Cal C(\bigtriangledown_2)$:
$\Cal C(\bigtriangledown_2)\simeq\overline{\Cal
C(\bigtriangledown_1)}\simeq\Cal C(\bigtriangledown_1)^{-1}
\simeq\,^{-1}\Cal C(\bigtriangledown_1)$,
$\Cal C(\bigtriangledown_1)\simeq\overline{\Cal
C(\bigtriangledown_2)}\simeq\Cal C(\bigtriangledown_2)^{-1}
\simeq\,^{-1}\Cal C(\bigtriangledown_2)$.


}
\enddocument